\documentclass[12pt,letterpaper]{amsart}
\usepackage{pslatex}
\usepackage[T1]{fontenc}
\usepackage[latin1]{inputenc}
\usepackage{amssymb}

\makeatletter
 \theoremstyle{plain}    
 \newtheorem{thm}{Theorem}[section]
 \numberwithin{equation}{section} %% Comment out for sequentially-numbered
 \numberwithin{figure}{section} %% Comment out for sequentially-numbered
 \theoremstyle{plain}
 \theoremstyle{remark}
 \newtheorem{rem}[thm]{Remark}
 \theoremstyle{definition}
 \newtheorem{defn}[thm]{Definition}
 \theoremstyle{remark}
 \newtheorem*{rem*}{Remark}
 \theoremstyle{plain}    
 \newtheorem*{prop*}{Proposition} 
 \theoremstyle{definition}
  \newtheorem{example}[thm]{Example}
 \theoremstyle{plain}    
 \newtheorem{prop}[thm]{Proposition} %%Delete [thm] to re-start numbering
 \theoremstyle{plain}    
 \newtheorem{lem}[thm]{Lemma} %%Delete [thm] to re-start numbering
 \theoremstyle{definition}
  \newtheorem{condition}[thm]{Condition}
 \theoremstyle{plain}    
 \newtheorem{cor}[thm]{Corollary} %%Delete [thm] to re-start numbering

\usepackage{verbatim} 
\usepackage{varioref}
\usepackage{eucal}
\usepackage{amsthm}
\usepackage[all]{xy}
\usepackage{rotating}
\usepackage{psfrag}
\usepackage{float}
\usepackage{dsfont}
\mathcode`\:="603A % Redefine ':' to be math punctuation since I never use
\DeclareSymbolFont{rsfs}{U}{rsfs}{m}{n}
\DeclareSymbolFontAlphabet{\mathrf}{rsfs}

\newcommand{\im}{\operatorname{im}}

\SelectTips{cm}{}
\xyoption{line}

\usepackage[english,american]{babel}
\makeatother
\begin{document}

\title{Cofree coalgebras over operads II.\\
Homology invariance}

\author{Justin R. Smith}

\subjclass{18D50;  Secondary: 16W30}

\keywords{operads, cofree coalgebras, homology invariance}

\address{Department of Mathematics\\
Drexel University\\
Philadelphia,~PA 19104}

\email{jsmith@drexel.edu}

\urladdr{http://vorpal.mcs.drexel.edu}

\date{\today}

\begin{abstract}
This paper gives conditions under which the cofree coalgebras constructed
in \cite{Smith-cofree} are homology invariant. 
\end{abstract}
\maketitle
\tableofcontents{}

\newcommand{\ring}{R}

\newcommand{\integers}{\mathbb{Z}}

\newcommand{\betabar}{\bar{\beta}}
 
\newcommand{\desusp}{\downarrow}

\newcommand{\susp}{\uparrow}

\newcommand{\cobar}{\mathcal{F}}

\newcommand{\bigp}{\mathop{\prod^{\prime}}}

\newcommand{\mfrac}{\mathfrak{M}}

\newcommand{\coend}{\mathrm{CoEnd}}

\newcommand{\ainfty}{A_{\infty}}

\newcommand{\coassoc}{\mathrm{Coassoc}}

\newcommand{\trm}{\mathrm{T}}

\newcommand{\tfr}{\mathfrak{T}}

\newcommand{\tabbr}{\hat{\trm}}

\newcommand{\Tabbr}{\hat{\tfr}}

\newcommand{\afr}{\mathfrak{A}}

\newcommand{\homz}{\mathrm{Hom}_{\ring}}

\newcommand{\homa}{\mathrm{Hom}}

\newcommand{\zend}{\mathrm{End}}

\newcommand{\rs}[1]{\mathrm{R}S_{#1 }}

\newcommand{\highprod}[1]{\bar{\mu}_{#1 }}

\newcommand{\slength}[1]{|#1 |}

\newcommand{\barcs}{\bar{\mathcal{B}}}

\newcommand{\ubarcs}{\mathcal{B}}

\newcommand{\zs}[1]{\mathbb{Z}S_{#1 }}

\newcommand{\homzs}[1]{\mathrm{Hom}_{\ring S_{#1 }}}

\newcommand{\zpi}{\mathbb{Z}\pi}

\newcommand{\D}{\mathfrak{D}}

\newcommand{\ahat}{\hat{\mathfrak{A}}}

\newcommand{\cbar}{{\bar{C}}}

\newcommand{\cf}[1]{\mathcal{C}(#1 )}

\newcommand{\ddelta}{\dot{\Delta}}

\newcommand{\dimlimiter}{\triangleright}

\newcommand{\coalgcat}{\mathrf S_{0}}

\newcommand{\hcoalgcat}{\mathrf{S}}

\newcommand{\ircoalgcat}{\mathrf I_{0}}

\newcommand{\bircoalgcat}{\mathrf{I}_{0}^{+}}

\newcommand{\hircoalgcat}{\mathrf I}

\newcommand{\chaincat}{\mathbf{Ch}(\ring)}

\newcommand{\coll}{\mathrm{Coll}}

\newcommand{\bchaincat}{\mathbf{Ch}^{+}(\ring)}

\newcommand{\ilimit}{\varprojlim\,}

\newcommand{\bigboxtimes}{\mathop{\boxtimes}}

\newcommand{\dlimit}{\varinjlim\,}

\newcommand{\coker}{\mathrm{{coker}}}

\newcommand{\iircoalgcat}{\mathrm{inv-}\ircoalgcat}

\newcommand{\dircoalgcat}{\mathrm{dir-}\ircoalgcat}

\newcommand{\core}[1]{\left\langle #1\right\rangle }

\newcommand{\ilimitder}{\varprojlim^{1}\,}

\newcommand{\homotopycat}{\mathbf{K}}

\newcommand{\homotopcellular}{\homotopycat_{\mathrm{cell}}}

\newcommand{\homotophom}{\mathrm{hom}_{\homotopycat(\integers)}}

\newcommand{\surj}{\mathrm{Surj}}

\newcommand{\surjf}{\surj_{f}}

\newcommand{\surjff}{\surjf^{2}}

\newcommand{\sethom}{\underline{\mathrm{Hom}}}

\newcommand{\setsethom}{\underline{\underline{\mathrm{Hom}}}}

\newcommand{\sethomn}[1]{\sethom_{#1}}

\newcommand{\setf}{\mathrm{Set}_{f}}

\newcommand{\setmod}{\setf\mathrm{-mod}}

\newcommand{\sigmamod}{\Sigma\mathrm{-mod}}

\newcommand{\order}{\mathrm{Ord}}

\newcommand{\coequalizer}{\mathop{\mathrm{coequalizer}}}

\newcommand{\unitobj}{\mathds{1}}

\newcommand{\forgetful}[1]{\lceil#1\rceil}

\newcommand{\pcoalg}[2]{P_{\mathcal{#1}}(#2) }

\newcommand{\mainoperad}{\mathcal{H}}

\section{Introduction}

The paper \cite{Smith-cofree} constructed cofree coalgebras over
operads cogenerated by free chain-complexes over a ring $\ring$.
The underlying chain-complexes of these cofree coalgebras were not
known to be free in the case where $\ring=\integers$ since they were
only submodules of the Baer-Specker group, $\integers^{\aleph_{0}}$
--- see \cite{Baer-Specker-survey} for a survey of this group.

In the present paper we address several issues:

\begin{enumerate}
\item We extend the construction of cofree coalgebras to the class of nearly
free modules --- see definition~\ref{def:nearlyfree} and appendix~\ref{sec:nearlyfree}.
This class includes free modules but is closed under the operations
of taking countable products and cofree coalgebras. Consequently,
it will be possible to iterate our cofree coalgebra construction.
\item We show that, under fairly weak conditions on the operad --- that
it is composed of projective modules that are finitely generated in
each dimension --- cofree coalgebras of nearly free chain-complexes
are homology invariant.
\end{enumerate}
Section~\ref{sec:generalconstruction} defines nearly free modules
and other terms connected with operads and coalgebras over them.

Section~\ref{sec:coalgnearlyfree} carries out step 1 above. It essentially
shows that cofree coalgebras preserve direct limits. Since nearly
free modules are direct limits of free modules, this defines cofree
coalgebras over nearly free modules.

Section~\ref{sec:Cofibrant-Operads} shows that cofibrant operads
are homotopy functors --- i.e. homotopies of maps induce homotopies
of cofree coalgebra morphisms. This, coupled with the results of appendix~\ref{sec:homotopyanddirect}
implies that they preserve homology equivalences of nearly free chain-complexes.

Our main result, proved in section~\ref{sec:The-general-case} is:

\medskip

\emph{Corollary~\ref{cor:mainresult}: Let $\ring$ be a field or
$\integers$ and let $\mathcal{V}=\{\mathcal{V}(n)\}$ be an operad
such that $\mathcal{V}(n)$ is $\ring S_{n}$-projective and finitely
generated in each dimension for all $n>0$. If}{\emph{\[
W_{\mathcal{V}}C=\left\{ \begin{array}{c}
L_{\mathcal{V}}C\\
M_{\mathcal{V}}C\\
P_{\mathcal{V}}C\\
\mathrf{F}_{\mathcal{V}}C\end{array}\right\} \]
}} \emph{--- the cofree coalgebras defined in \cite{Smith-cofree}
--- and\[
f:C\to D\]
is a homology equivalence of nearly free chain complexes (see definition~\ref{def:nearlyfree})
that are bounded from below, then the induced map\[
W_{\mathcal{V}}f:W_{\mathcal{V}}C\to W_{\mathcal{V}}D\]
is a homology equivalence. }

\begin{rem}
The condition on $\mathcal{V}$ is essentially equivalent to the condition
of being $\Sigma$-cofibrant in \cite{Berger-moerdijk-axiom-operad}. 

This condition is necessary because there are well-known cases in
which it does not hold and the associated cofree coalgebras are \emph{not}
homology invariant.
\end{rem}
\newdir{ >}{{}*!/-5pt/@{>}}

\section{\label{sec:generalconstruction}Definitions}

Throughout this paper, $\ring$ will denote a field or $\integers$. 

\begin{defn}
\label{def:nearlyfree}An $\ring$-module $M$ will be called \emph{nearly
free} if every countable submodule is $\ring$-free.
\end{defn}
\begin{rem*}
This condition is automatically satisfied unless $\ring=\integers$.

Clearly, any $\integers$-free module is also nearly free. The Baer-Specker
group, $\integers^{\aleph_{0}}$, is a well-known example of a nearly
free $\integers$-module that is \emph{not} free --- see \cite{Fuchs:1970},
\cite{Baer:1937}, and \cite{Baer-Specker-nonfree}. Compare this
with the notion of $\aleph_{1}$-\emph{free} \emph{groups} --- see
\cite{Blass-Gobel:1996}. 

By abuse of notation, we will often call chain-complexes nearly free
if their underlying modules are (when one ignores grading).

Nearly free $\integers$-modules enjoy useful properties that free
modules do \emph{not}. For instance, in many interesting cases, the
cofree coalgebra of a nearly free chain-complex is nearly free.
\end{rem*}
We will denote the closed symmetric monoidal category of (not necessarily
free) $\ring$-chain-complexes with $\ring$-tensor products by $\chaincat$.
These chain-complexes are allowed to extend into arbitrarily many
negative dimensions and have underlying graded $\ring$-modules that
are

\begin{itemize}
\item arbitrary if $\ring$ is a field (but they will be free)
\item \emph{nearly free,} in the sense of definition~\ref{def:nearlyfree},
if $\ring=\integers$. 
\end{itemize}
\begin{defn}
\label{def:unitinterval}The object $I\in\chaincat$, the \emph{unit
interval,} is defined by\[
I_{k}=\begin{cases}
\ring\cdot p_{0}\oplus R\cdot p_{1} & \mathrm{if\,}k=0\\
\ring\cdot q & \mathrm{if\,}k=1\\
0 & \mathrm{if\,}k\ne0,1\end{cases}\]
where $p_{0}$, $p_{1}$, $q$ are just names for the canonical generators
of $I$, and the one nonzero boundary map is defined by $q\mapsto p_{1}-p_{0}$.

We also define, for any object $A\in\chaincat$, the \emph{cone on}
$A$, denoted $\bar{A}$ and equal to $A\otimes I/A\otimes R\cdot p_{1}$.
There are canonical morphisms $A\to\bar{A}$ and $\bar{A}\to\Sigma A$,
where $\Sigma:\chaincat\to\chaincat$ is the functor that raises the
grading by 1.
\end{defn}
Two morphisms\[
f_{0},\, f_{1}:C\to D\]
in $\chaincat$, are defined to be \emph{chain-homotopic} if there
exists a morphism\[
F:C\otimes I\to D\]
such that $F|C\otimes R\cdot p_{i}=f_{i}:C\to D$. This is well-known
to be equivalent to the existence of a degree $+1$ map $\Phi:C\to D$
such that $\partial_{D}\circ\Phi+\Phi\circ\partial_{C}=f_{1}-f_{0}$.

We make extensive use of the Koszul Convention (see~\cite{Gugenheim:1960})
regarding signs in homological calculations:

\begin{defn}
\label{def:koszul} If $f:C_{1}\to D_{1}$, $g:C_{2}\to D_{2}$ are
maps, and $a\otimes b\in C_{1}\otimes C_{2}$ (where $a$ is a homogeneous
element), then $(f\otimes g)(a\otimes b)$ is defined to be $(-1)^{\deg(g)\cdot\deg(a)}f(a)\otimes g(b)$. 
\end{defn}
\begin{rem}
If $f_{i}$, $g_{i}$ are maps, it isn't hard to verify that the Koszul
convention implies that $(f_{1}\otimes g_{1})\circ(f_{2}\otimes g_{2})=(-1)^{\deg(f_{2})\cdot\deg(g_{1})}(f_{1}\circ f_{2}\otimes g_{1}\circ g_{2})$.
\end{rem}
\begin{defn}
\label{def:homcomplex}Given chain-complexes $A,B\in\chaincat$ define\[
\homz(A,B)\]
 to be the chain-complex of graded $\ring$-morphisms where the degree
of an element $x\in\homz(A,B)$ is its degree as a map and with differential\[
\partial f=f\circ\partial_{A}-(-1)^{\deg f}\partial_{B}\circ f\]
As a $\ring$-module $\homz(A,B)_{k}=\prod_{j}\homz(A_{j},B_{j+k})$.
\end{defn}
\begin{rem*}
Given $A,B\in\chaincat^{S_{n}}$, we can define $\homzs{n}(A,B)$
in a corresponding way.
\end{rem*}

\begin{defn}
\label{def:setfmodules}Define:

\begin{enumerate}
\item $\setf$ to be the category of finite sets and bijections. Let $\setf^{2}$
be the category of finite sets whose elements are also finite sets.
Morphisms are bijections of sets that respect the {}``fine structure''
of elements that are also sets. There is a \emph{forgetful functor}\[
\mathfrak{f}:\setf^{2}\to\setf\]
 that simply forgets that the elements of an object of $\setf^{2}$
are, themselves, finite sets. There is also a \emph{{}``flattening''
functor}\[
\mathfrak{g}:\setf^{2}\to\setf\]
 that sends a set (of sets) to the union of the elements (regarded
as sets).
\item For a finite set $X$, $\Sigma_{X}=\mathrm{End}_{\setf}(X)$.
\item $\setmod$ to be the category of contravariant functors $\mathrm{Func}(\setf^{\mathrm{op}},\chaincat)$,
with morphisms that are natural transformations. 
\item Given $C,\, D\in\setmod$, define $\sethom(C,D)$ to be the set of
natural transformations of functors. Also define $\sethomn{X}(C,D)$,
where $X\in\setf$, to be the natural transformations of $C$ and
$D$ restricted to sets isomorphic to $X$ (i.e., of the same cardinality).
Both of these functors are chain-complexes.
\item $\sigmamod$ to be the category of sequences $\{ M(n)\}$, $m\ge1$
where $M(n)\in\chaincat$ and $M(n)$ is equipped with a right $S_{n}$-action.
\end{enumerate}
\end{defn}
\begin{rem*}
If $[n]$ is the set of the first $n$ positive integers, then $\Sigma_{[n]}=S_{n}$,
the symmetric group. If $M$ is a $\setf$-module then, for each finite
set, $X$, there is a right $\Sigma_{X}$-action on $M(X)$.

We follow the convention that $S_{0}=S_{1}=\{1\}$, the trivial group.

Note that $\sigmamod$ is what is often called the category of collections.

If $\mathbf{a}=\{\{ x\},\{ y,z,t\},\{ h\}\}\in\setf^{2}$ then $\mathfrak{f}(\mathbf{a})\cong[3]$,
a set of three elements, and $\mathfrak{g}(\mathbf{a})=\{ x,y,z,t,h\}$.
\end{rem*}
It is well-known that the categories $\setmod$ and $\sigmamod$ are
isomorphic --- see section~1.7 in part~I of \cite{operad-book}.
The restriction isomorphism\[
r:\setmod\to\sigmamod\]
simply involves evaluating functors on the finite sets $[n]$ for
all $n\ge1$. If $F\in\setmod$, then $r(F)=\{ F([n])\}$. The functorial
nature of $F$ implies that $F([n])$ is equipped with a natural $S_{n}$-action.
The functors $\sethomn{n}(C,D)$ correspond to $\homzs{n}(C([n]),D([n]))$
and the fact that morphisms in $\setf$ preserve cardinality imply
that\[
\sethom(C,D)=\prod_{n\ge0}\sethomn{n}(C,D)\]

Although $\setf$-modules are equivalent to modules with a symmetric
group action, it is often easier to formulate operadic constructions
in terms of $\setmod$. Equivariance relations are automatically satisfied.

\begin{defn}
\label{def:orderings}If $X$ is a finite set of cardinality $n$
the \emph{set of orderings of} $X$ is \[
\order(X)=\{ f|f:X\xrightarrow{\cong}[n]\}\]

\end{defn}
Now we define a $\setf$ analogue to the multiple tensor product.
Given a set $X$ of cardinality $n$, and an assignment of an object
$C_{x}\in\chaincat$ for each element $x\in X$, we can define, for
each $g\in\order(X)$ a product\[
\bigotimes_{g}C_{x}=C_{g^{-1}(1)}\otimes\cdots\otimes C_{g^{-1}(n)}\]
The symmetry of tensor products determines a morphism \[
\bar{\sigma:}\bigotimes_{g}C_{x}\to\bigotimes_{\sigma\circ g}C_{x}\]
for each $\sigma\in S_{n}$ which essentially permutes factors and
multiplies by $\pm1$, following the Koszul Convention in definition~\ref{def:koszul}.

\begin{defn}
\label{def:unorderedtensor}The \emph{unordered tensor product} is
defined by\[
\bigotimes_{X}C_{x}=\coequalizer_{\sigma\in S_{n}}\left\{ \bar{\sigma}:\bigoplus_{g\in\order(X)}\bigotimes_{g}C_{x}\to\bigoplus_{g\in\order(X)}\bigotimes_{g}C_{x}\right\} \]

If $C\in\chaincat$ and $X\in\setf$ then $C^{X}$ will denote the
unordered tensor product \[
\bigotimes_{X}C\]
of copies of $C$ indexed by elements of $X$, and $C^{\otimes}$will
denote the $\setf$-module whose value on $X\in\setf$ is $C^{X}$.

We use $X\cdot C$ to denote a direct sum of $n$ copies of $C$,
where $n$ is the cardinality of a finite set $X$.

When $X\in\setf^{2}$, \[
\bigotimes_{X}C\]
 is regarded as being taken over $\mathfrak{f}(X)$ --- i.e., we {}``forget''
that the elements of $X$ are sets themselves.
\end{defn}
\begin{rem*}
The unordered tensor product is isomorphic (as an object of $\chaincat$
to the tensor product of the $C_{x}$, as $x$ runs over the elements
of $X$. The coequalizer construction determines how the it behaves
with respect to \emph{set-morphisms.}

If $X=[n]$, then $C^{[n]}=C^{n}$. Note that $C^{X}\otimes C^{Y}=C^{X\sqcup Y}$,
for $X,\, Y\in\setf$. We also follow the convention that $C^{\emptyset}=\mathds{1}=\ring$,
concentrated in dimension $0$.
\end{rem*}
\begin{defn}
\label{def:unorderedspans}If $X\in\setf$, $x\in X$ and $\{ f_{y}:V_{y}\to U_{y}\}$are
morphisms of $\chaincat$ indexed by elements $y\in X$ then define\[
\bigotimes_{X,x}(U,V)=\bigotimes_{y\in X}Z_{y}\xrightarrow{1\otimes\cdots\otimes f_{x}\otimes\cdots\otimes1}\bigotimes_{y\in X}U_{y}\]
 to be the unordered tensor product, where \[
Z_{y}=\begin{cases}
U_{y} & \text{if\,}y\ne x\\
V_{y} & \text{if\,}u=x\end{cases}\]

\end{defn}
\begin{rem*}
Given any ordering of the elements of the set $X$, there exists a
canonical isomorphism\[
\bigotimes_{X,x}(U,V)=\underbrace{U\otimes\cdots\otimes V\otimes\cdots\otimes U}_{\text{position\,}x}\]

\end{rem*}
\begin{defn}
\label{def:setxunion}Let $X,Y\in\setf$ and let $x\in X$. Define\[
X\sqcup_{x}Y=(X\setminus\{ x\})\cup Y\]

\end{defn}
\begin{rem*}
Note that $X\sqcup_{x}\emptyset=X\setminus\{ x\}$.
\end{rem*}
\begin{prop*}
\label{pro:setxunionprops}If $X,\, Y,\, Z\in\setf$, and $x\, x_{1},\, X_{2}\in X$
and $y\in Y$, then\begin{eqnarray*}
X\sqcup_{x}(Y\sqcup_{y}Z) & = & (X\sqcup_{x}Y)\sqcup_{y}Z\\
(X\sqcup_{x_{1}}Y)\sqcup_{x_{2}}Z & = & (X\sqcup_{x_{2}}Z)\sqcup_{x_{1}}Y\end{eqnarray*}

\end{prop*}
\begin{defn}
\label{def:operad}An \emph{operad} in $\chaincat$ is a $\setf$-module,
$C$ equipped with operations\[
\circ_{x}:C(X)\otimes C(Y)\to C(X\sqcup_{x}Y)\]
for all $x\in X$and all $X,\, Y\in\setf$ and satisfying the two
axioms
\begin{enumerate}
\item \emph{Associativity:}\begin{multline*}
\circ_{x}(1\otimes\circ_{y})=\circ_{y}(\circ_{x}\otimes1):\\
C(X)\otimes C(Y)\otimes C(Z)\to C(X\sqcup_{x}(Y\sqcup_{y}Z))\end{multline*}
\begin{multline*}
\circ_{x_{2}}(\circ_{x_{1}}\otimes1)=\circ_{x_{1}}(\circ_{x_{2}}\otimes1)(1\otimes\tau):\\
C(X)\otimes C(Y)\otimes C(Z)\to C((X\sqcup_{x_{1}}Y)\sqcup_{x_{2}}Z)\end{multline*}
for all $X,\, Y,\, Z\in\setf$ and all $x\, x_{1},\, x_{2}\in X$
and $y\in Y$, where $\tau:C(Y)\otimes C(Z)\to C(Z)\otimes C(Y)$
is the transposition isomorphism.
\item \emph{Unit:} There exist morphisms $\eta_{x}:\mathds{1}\to C(\{ x\})$
for all singleton sets $\{ x\}\in\setf$ that make the diagrams\[\xymatrix{{C(X)\otimes\mathds{1}}\ar[r]^{\cong}\ar[d]_{1\otimes\eta_x} & {C(X)} \\ {C(X)\otimes C({x})}\ar[ru]_{\circ_x} & {}} \qquad\xymatrix{{\mathds{1}\otimes C(X)}\ar[r]^{\cong}\ar[d]_{\eta_x\otimes 1} & {C(X)} \\ {C(X)}\ar[ru]_{\circ_x} & {}}\]commute,
for all $X\in\setf$. The operad will be called \emph{nonunital} if
the axioms above only hold for \emph{nonempty} sets.
\end{enumerate}
\end{defn}
\begin{rem*}
See theorem~1.60 and 1.61 and section 1.7.1 of \cite{operad-book}
for the proof that this defines operads correctly. For more traditional
definitions, see \cite{Smith-cofree}, \cite{Kriz-May}. This is basically
the definition of a pseudo-operad in \cite{operad-book} where we
have added the unit axiom. To translate this definition into the more
traditional ones, set the $n^{\mathrm{th}}$ \emph{component} of the
operad to $C([n])$.

The use of $\setmod$ causes the equivariance conditions in \cite{Kriz-May}
to be automatically satisfied. 

The operads we consider here correspond to \emph{symmetric} operads
in \cite{Smith-cofree}.

The term {}``unital operad'' is used in different ways by different
authors. We use it in the sense of Kriz and May in \cite{Kriz-May},
meaning the operad has a $0$-component that acts like an arity-lowering
augmentation under compositions. This is $C(\emptyset)=\mathds{1}$.
\end{rem*}
A simple example of an operad is:

\begin{example}
\label{example:frakS0}For each finite set, $X$ $,C(X)=\integers\Sigma_{X}$,
with composition defined by inclusion of sets. This operad is denoted
$\mathfrak{S}_{0}$. In other notation, its $n^{\text{th}}$component
is the \emph{symmetric group-ring} $\integers S_{n}$.
\end{example}
For the purposes of this paper, the canonical example of an operad
is

\begin{defn}
\label{def:coend}Given any $C\in\chaincat$, the associated \emph{coendomorphism
operad}, $\coend(C)$ is defined by\[
\coend(C)(X)=\homz(C,C^{X})\]
for $X\in\setf$, and $C^{X}=\bigotimes_{X}C$ is the unordered tensor
product defined in definition~\ref{def:unorderedtensor}. The compositions
$\{\circ_{x}\}$ are defined by\begin{multline*}
\circ_{x}:\homz(C,C^{X})\otimes\homz(C,C^{Y})\to\\
\homz(C,C^{X\setminus\{ x\}}\otimes C_{x}\otimes\homz(C,C^{Y}))\xrightarrow{\homz(1,1\otimes e)}\\
\homz(C,C^{X\setminus\{ x\}}\otimes C^{Y}))=\homz(C,C^{X\sqcup_{x}Y})\end{multline*}
where $C_{x}$ is the copy of $C$ corresponding to $x\in X$ and
$e:C_{x}\otimes\homz(C,C^{Y})\to C^{Y}$ is the evaluation morphism.
This is a non-unital operad, but if $C\in\chaincat$ has an augmentation
map $\varepsilon:C\to\mathds{1}$ then we can set\[
\coend(C)(\emptyset)=\mathds{1}\]
and \begin{multline*}
\circ_{x}:\homz(C,C^{X})\otimes\homz(C,C^{\emptyset})=\homz(C,C^{X})\otimes\mathds{1}\\
\xrightarrow{\homz(1,1_{X\setminus\{ x\}}\otimes\varepsilon_{x})}\homz(C,C^{X\setminus\{ x\}})\end{multline*}
where $1_{X\setminus\{ x\}}:C^{X\setminus\{ x\}}\to C^{X\setminus\{ x\}}$
is the identity map and $\varepsilon_{x}:C_{x}\to\mathds{1}$ is the
augmentation, applied to the copy of $C$ indexed by $x\in X$.

Given $C\in\chaincat$ with \emph{subcomplexes} $\{ D_{1},\dots,D_{k}\}$,
the \emph{relative coendomorphism operad} $\coend(C;\{ D_{i}\})$
is defined to be the sub-operad of $\coend(C)$ consisting of maps
$f\in\homz(C,C^{X})$ such that $f(D_{j})\subseteq D_{j}^{X}\subseteq C^{X}$
for all $j$.
\end{defn}
We use the coendomorphism operad to define the main object of this
paper:

\begin{defn}
\label{def:coalg}A \emph{coalgebra over an operad} $\mathcal{V}$
is a chain-complex $C\in\chaincat$ with an operad morphism $\alpha:\mathcal{V}\to\coend(C)$,
called its \emph{structure map.} We will sometimes want to define
coalgebras using the \emph{adjoint structure map}\[
\bar{\alpha}:C\to\sethom(\mathcal{V},C^{\otimes})\]
(in $\chaincat)$ or even the set of chain-maps\[
\bar{\alpha}_{X}:C\to\sethomn{X}(\mathcal{V}(X),C^{X})\]
for all $X\in\setf$.
\end{defn}
We can also define the analogue of an ideal:

\begin{defn}
\label{def:coideal}Let $C$ be a coalgebra over the operad $\mathcal{U}$
with adjoint structure map \[
\alpha:C\to\sethom(\mathcal{U},C^{\otimes})\]
and let $D\subseteq\forgetful{C}$ be a sub-chain complex that is
a direct summand. Then $D$ will be called a \emph{coideal} of $C$
if the composite \[
\alpha|D:D\to\sethom(\mathcal{U},C^{\otimes})\xrightarrow{\sethom(1_{\mathcal{U}},p^{\otimes})}\sethom(\mathcal{U},(C/D)^{\otimes})\]
 \emph{vanishes}, where $p:C\to C/D$ is the projection to the quotient
(in $\chaincat$). 
\end{defn}
\begin{rem*}
Note that it is easier for a sub-chain-complex to be a coideal of
a coalgebra than to be an ideal of an algebra. For instance, all sub-coalgebras
of a coalgebra are also coideals. Consequently it is easy to form
quotients of coalgebras and hard to form sub-coalgebras. This is dual
to what occurs for algebras.
\end{rem*}

We will sometimes want to focus on a particular class of $\mathcal{V}$-coalgebras:
the \emph{pointed, irreducible coalgebras}. We define this concept
in a way that extends the conventional definition in \cite{Sweedler:1969}:

\begin{defn}
\label{def:pointedirreducible} Given a coalgebra over a unital operad
$\mathcal{V}$ with adjoint structure-map\[
a_{X}:C\to\sethomn{X}(\mathcal{V}(X),C^{X})\]
an element $c\in C$ is called \emph{group-like} if $a_{X}(c)=f_{X}(c^{X})$
for all $n>0$. Here $c^{X}\in C^{X}$ is the $n$-fold $\ring$-tensor
product, where $n$ is the cardinality of $X$,\[
f_{X}=\homz(\epsilon_{X},1):\homz(\mathds{1},C^{X})=C^{X}\to\sethomn{X}(\mathcal{V}(X),C^{X})\]
and $\epsilon_{X}:\mathcal{V}(X)\to\mathcal{V}(\emptyset)=\mathds{1}=\ring$
is the augmentation (which is $n$-fold composition with $\mathcal{V}(\emptyset)$). 

A coalgebra $C$ over an operad $\mathcal{V}$ is called \emph{pointed}
if it has a \emph{unique} group-like element (denoted $1$), and \emph{pointed
irreducible} if the intersection of any two sub-coalgebras contains
this unique group-like element.
\end{defn}
\begin{rem*}
Note that a group-like element generates a sub $\mathcal{V}$-coalgebra
of $C$ and must lie in dimension $0$.

Although this definition seems contrived, it arises in {}``nature'':
The chain-complex of a pointed, simply-connected reduced simplicial
set is naturally a pointed irreducible coalgebra over the Barratt-Eccles
operad, $\mathfrak{S}=\{ C(K(S_{n},1))\}$ (see \cite{Smith:1994}).
In this case, the operad action encodes the chain-level effect of
Steenrod operations.
\end{rem*}
\begin{prop}
Let $D$ be a pointed, irreducible coalgebra over an operad $\mathcal{V}$.
Then the augmentation map\[
\varepsilon:D\to\ring\]
is naturally split and any morphism of pointed, irreducible coalgebras
\[
f:D_{1}\to D_{2}\]
 is of the form\[
1\oplus\bar{f}:D_{1}=\ring\oplus\ker\varepsilon_{D_{1}}\to D_{2}=\ring\oplus\ker\varepsilon_{D_{2}}\]
where $\varepsilon_{i}:D_{i}\to\ring$, $i=1,2$ are the augmentations.
\end{prop}
\begin{proof}
The definition (\ref{def:pointedirreducible}) of the sub-coalgebra
$\ring\cdot1\subseteq D_{i}$ is stated in an invariant way, so that
any coalgebra morphism must preserve it. Any morphism must also preserve
augmentations because the augmentation is the $0^{\mathrm{th}}$-order
structure-map. Consequently, $f$ must map $\ker\varepsilon_{D_{1}}$to
$\ker\varepsilon_{D_{2}}$. The conclusion follows.
\end{proof}
\begin{defn}
\label{def:pointedirredcat} We denote the \emph{category} of coalgebras
over $\mathcal{V}$ by $\coalgcat$. If $\mathcal{V}$ is unital,
every $\mathcal{V}$-coalgebra, $C$, comes equipped with a canonical
augmentation\[
\varepsilon:C\to\ring\]
so the \emph{terminal object} is $\ring$. If $\mathcal{V}$ is not
unital, the terminal object \emph{}in this category is $0$, the null
coalgebra.

The category of \emph{pointed irreducible coalgebras} over $\mathcal{V}$
is denoted $\ircoalgcat$ --- this is only defined if $\mathcal{V}$
is unital. Its terminal object is the coalgebra whose underlying chain
complex is $\ring$ concentrated in dimension $0$.
\end{defn}
We also need:

\begin{defn}
\label{def:forgetful}If $A\in\mathrf C=\ircoalgcat$ or $\coalgcat$,
then $\forgetful{A}$ denotes the underlying chain-complex in $\chaincat$
of\[
\ker A\to t\]
where $t$ denotes the terminal object in $\mathrf{C}$ --- see definition~
\ref{def:pointedirredcat}. We will call $\forgetful{\ast}$ the \emph{forgetful
functor} from $\mathrf{C}$ to $\chaincat$.
\end{defn}
We will use the concept of cofree coalgebra cogenerated by a chain
complex: 

\begin{defn}
\label{def:cofreecoalgebra}Let $D$ be a coalgebra over an operad
$\mathcal{U}$, equipped with a $\chaincat$-morphism $\varepsilon:\forgetful{D}\to E$,
where $E\in\chaincat$. Then $D$ is called \emph{the cofree coalgebra
over} $\mathcal{U}$ \emph{cogenerated} \emph{by} $\varepsilon$ if
any morphism in $\chaincat$\[
f:\forgetful{C}\to E\]
where $C$ is a $\mathcal{U}$-coalgebra, induces a \emph{unique}
morphism of $\mathcal{U}$-coalgebras\[
\alpha_{f}:C\to D\]
that makes the diagram \[\xymatrix{{\forgetful{C}}\ar[r]^{\forgetful{\alpha_f}}\ar[rd]_{f} & {\forgetful{D}}\ar[d]^{\varepsilon}\\ {} & {E}}\]
commute. Here $\alpha_{f}$ is called the \emph{classifying map} of
$f$. If $C$ is a $\mathcal{U}$-coalgebra then \[
\alpha_{1}:C\to L_{\mathcal{U}}\forgetful{C}\]
 will be called the \emph{classifying map of} $C$.
\end{defn}
This universal property of cofree coalgebras implies that they are
unique up to  isomorphism if they exist.

\section{Extending the construction in \cite{Smith-cofree}\label{sec:coalgnearlyfree}}

The paper \cite{Smith-cofree} gave an explicit construction of $L_{\mathcal{U}}C$
when $C$ was an $\ring$-free chain complex. When $\ring$ is a field,
all chain-complexes are $\ring$-free, so the results of the present
paper are already true in that case. 

Consequently, we will restrict ourselves to the case where $\ring=\integers$. 

\begin{prop}
The forgetful functor (defined in definition~\ref{def:forgetful})
and cofree coalgebra functors define adjoint pairs\begin{eqnarray*}
\pcoalg{V}{\ast}:\chaincat & \leftrightarrows & \ircoalgcat:\forgetful{\ast}\\
L_{\mathcal{V}}(\ast):\chaincat & \leftrightarrows & \coalgcat:\forgetful{\ast}\end{eqnarray*}

\end{prop}
\begin{rem*}
The adjointness of the functors follows from the universal property
of cofree coalgebras --- see~\cite{Smith-cofree}. 
\end{rem*}
The Adjoints and Limits Theorem in \cite{MacLane:cw} implies that:

\begin{thm}
\label{th:cofreelimits}If $\{ A_{i}\}$ is an inverse system in $\chaincat$
and $\{ C_{i}\}$ is a direct system in $\ircoalgcat$ or $\coalgcat$
then\begin{eqnarray*}
\ilimit P_{\mathcal{V}}(A_{i}) & = & P_{\mathcal{V}}(\ilimit A_{i})\\
\ilimit L_{\mathcal{V}}(A_{i}) & = & L_{\mathcal{V}}(\ilimit A_{i})\\
\forgetful{\dlimit C_{i}} & = & \dlimit\forgetful{C_{i}}\end{eqnarray*}

\end{thm}
\begin{rem*}
This implies that \emph{direct} limits in $\ircoalgcat$ or $\coalgcat$
are the same as direct limits of underlying chain-complexes. 
\end{rem*}
\begin{prop}
\label{pro:directlimitfintelygenerated}If $C\in\chaincat$, let $\mathrf{G}(C)$
denote the lattice of countable subcomplexes of $C$. Then\[
C=\dlimit\mathrf{G}(C)\]

\end{prop}
\begin{proof}
Clearly $\dlimit\mathrf{G}(C)\subseteq C$ since all of the canonical
maps to $C$ are inclusions. Equality follows from every element $x\in C$
being contained in a finitely generated subcomplex of $C$ consisting
of $x$ and $\partial(x)$. 
\end{proof}
\begin{lem}
\label{lem:homdirectlimit}Let $n>1$ be an integers, let $F$ be
a finitely-generated projective (non-graded) $\zs{n}$-module, and
let $\{ C_{\alpha}\}$ a direct system of modules. Then the natural
map\[
\dlimit\homzs{n}(F,C_{\alpha})\to\homzs{n}(F,\dlimit C_{\alpha})\]
is an isomorphism.

If $F$ and the $\{ C_{\alpha}\}$ are graded, the corresponding statement
is true if $F$ is finitely-generated and $\zs{n}$-projective in
each dimension.
\end{lem}
\begin{proof}
We will only prove the non-graded case. The graded case follows from
the fact that the maps of the $\{ C_{\alpha}\}$ preserve grade.

In the non-graded case, finite generation of $F$ implies that the
natural map\[
\bigoplus_{\alpha}\homzs{n}(F,C_{\alpha})\to\homzs{n}(F,\bigoplus_{\alpha}C_{\alpha})\]
is an isomorphism. The projectivity of $F$ implies that $\homzs{n}(F,*)$
is exact, so the short exact sequence defining the direct limit is
preserved.
\end{proof}
\begin{thm}
\label{th:cofreedirectlimit}Let $\mathcal{V}=\{\mathcal{V}(X)\}$
be an operad and let $C$ be a chain-complex with $\mathrf{G}(C)=\{ C_{\alpha}\}$
the direct system of countable subcomplexes ordered by inclusion.
In addition, suppose:
\begin{enumerate}
\item For all $n\ge0$, $\mathcal{V}(X)$ is $\integers\Sigma_{X}$-projective
and finitely generated in each dimension.
\item $C$ is nearly free (see definition~\ref{def:nearlyfree}).
\end{enumerate}
Then the cofree coalgebras \[
L_{\mathcal{V}}C,\, P_{\mathcal{V}}C,\, M_{\mathcal{V}}C,\,\mathrf{F}_{\mathcal{V}}C\]
are well-defined and\[
\left.\begin{array}{c}
L_{\mathcal{V}}C=\dlimit L_{\mathcal{V}}C_{\alpha}\\
P_{\mathcal{V}}C=\dlimit P_{\mathcal{V}}C_{\alpha}\\
M_{\mathcal{V}}C=\dlimit M_{\mathcal{V}}C_{\alpha}\\
\mathrf{F}_{\mathcal{V}}C=\dlimit\mathrf{F}_{\mathcal{V}}C_{\alpha}\end{array}\right\} \subseteq\sethom(\mathcal{V},C^{\otimes})\]

\end{thm}
\begin{rem*}
Indeed, the construction of them given in \cite{Smith-cofree} is
valid in this case. 
\end{rem*}
\begin{proof}
The only part of the construction in \cite{Smith-cofree} that uses
$\integers$-freeness is the proof that the $L_{\mathcal{V}}C$ are
coalgebras --- i.e., that the diagrams in Appendix~B of \cite{Smith-cofree}
commute. The construction of the $L_{\mathcal{V}}C$ (as \emph{chain-complexes})
does not use it. 

The near-freeness of $C$ implies that the $C_{\alpha}$ are all free.

We will regard the \emph{chain-complex,} $\forgetful{L_{\mathcal{V}}C}$,
as the result of this construction in Lemma~3.4 of \cite{Smith-cofree}
--- setting aside questions of whether it's a coalgebra.

We have $C=\dlimit C_{\alpha}$ and the conditions on $\mathcal{V}$
(and lemma~\ref{lem:homdirectlimit}) imply that \[
\sethom(\mathcal{V},C^{\otimes})=\prod_{n>0}\homzs{n}(\mathcal{V}_{n},C^{n})=\prod_{n>0}\dlimit\homzs{n}(\mathcal{V}_{n},C_{\alpha_{n}}^{n})\]
where the $C_{\alpha_{n}}$ are countable. This and the $\integers$-flatness
of $C$ implies that every \[
x\in\prod_{n>0}\homzs{n}(\mathcal{V}_{n},C^{n})\]
lies in the image of \[
\prod_{n>0}\homzs{n}(\mathcal{V}_{n},C_{\alpha_{n}}^{n})\]
for some countable subcomplexes $\{ C_{\alpha_{n}}\}$.

We claim the natural map\[
\dlimit\forgetful{L_{\mathcal{V}}C_{\alpha}}\to\forgetful{L_{\mathcal{V}}C}\]
 is surjective. If $x\in\forgetful{L_{\mathcal{V}}C}$ is contained
in\[
\prod_{n>0}\homzs{n}(\mathcal{V}_{n},C_{\alpha_{n}}^{n})\subseteq\prod_{n>0}\homzs{n}(\mathcal{V}_{n},C^{n})\]
where the $\{ C_{\alpha_{n}}\}$ are all countable, then \[
\bar{C}=\sum_{n=1}^{\infty}C_{\alpha_{n}}\]
is also countable, and $x$ is in the image of an element $y\in\forgetful{L_{\mathcal{V}}\bar{C}}$.

Consequently\[
\forgetful{L_{\mathcal{V}}C}=\dlimit\forgetful{L_{\mathcal{V}}C_{\alpha}}\]
and theorem~\ref{th:cofreelimits} implies that this direct limit
has a natural coalgebra structure. The conclusion follows.
\end{proof}

\section{Cofibrant Operads\label{sec:Cofibrant-Operads}}

We define conditions on operads that ensure they are homotopy functors
and then apply the main result to show that they are homology invariant.

Now we determine the conditions necessary to make cofree coalgebras
into homotopy functors.

The relative coendomorphism operad of the unit interval is

\begin{condition}
\emph{\label{cond:mainoperadassumption}Throughout the rest of this
section,} we assume that $\mathcal{V}$ is an operad equipped with
a morphism of operads\[
\delta:\mathcal{V}\to\mathcal{V}\otimes\mathfrak{S}_{0}\]
(see definition~\ref{def:coend}) that makes the diagram \[\xymatrix{{\mathcal{V}}\ar[r]^-{\delta}\ar@{=}[rd] & {\mathcal{V}\otimes\mathfrak{S}_0}\ar[d] \\ {} & {\mathcal{V}}}\]
commute. Here, the operad structure on $\mathcal{V}\otimes\mathfrak{S}_{0}$
is just the tensor product of the operad structures of $\mathcal{V}$
and $\mathfrak{S}_{0}$. 

We also assume that the arity-$1$ component of $\mathcal{V}$ is
equal to $\ring$, generated by the unit.
\end{condition}
This is similar to the conditions satisfied by $\Sigma$-split operads
in \cite{Berger-moerdijk-axiom-operad}. It is also satisfied by cofibrant
operads --- of which the most straightforward class is that of free
operads.

The significance of this condition is given by:

\begin{prop}
\label{pro:ctimesi}Suppose the operad $\mathcal{V}$ satisfies Condition~\ref{cond:mainoperadassumption}
and $C$ is a $\mathcal{V}$-coalgebra. Then the coalgebra structure
of $C$ naturally extends to a coalgebra structure on $C\otimes I$
whose restrictions to $C\otimes p_{i}$, $i=0,1$ agree with the coalgebra
structure of $C$.
\end{prop}
\begin{proof}
The results of appendix~\ref{sec:relativecoendomorphism} imply that\[
\coend(I;\{\integers\cdot p_{0},\integers\cdot p_{1}\})=\mathfrak{S}_{0}\]
so that Condition~\ref{cond:mainoperadassumption} implies that the
operad morphism\[
\mathcal{V}\to\coend(C)\]
defining the coalgebra structure of $C$, lifts to a morphism\[
\mathcal{V}\to\coend(C)\otimes\coend(I;\{\integers\cdot p_{0},\integers\cdot p_{1}\})\to\coend(C\otimes I)\]
whose coalgebra-structure on $C\otimes\{ p_{i}\}$, $i=0,1$ coincides
with that of $C$.
\end{proof}
Our condition implies that:

\begin{prop}
\label{prop:lefthomotopy}Let $C$ and $D$ be objects of $\chaincat$
and let \[
f_{1},f_{2}:C\to D\]
 be chain-homotopic morphisms via a chain-homotopy \begin{equation}
F:C\otimes I\to D\label{eq:chainhomotop}\end{equation}
Then the induced maps\begin{eqnarray*}
P_{\mathcal{V}}f_{i}:P_{\mathcal{V}}C & \to & P_{\mathcal{V}}D\\
L_{\mathcal{V}}f_{i}:L_{\mathcal{V}}C & \to & L_{\mathcal{V}}D\end{eqnarray*}
$i=1,2$, are left-homotopic in $\ircoalgcat$ and $\coalgcat$, respectively
via a chain homotopy\[
F':P_{\mathcal{V}}f_{i}:(P_{\mathcal{V}}C)\otimes I\to P_{\mathcal{V}}D\]
If we equip $C\otimes I$ with a coalgebra structure using condition~\ref{cond:mainoperadassumption}
and proposition~\ref{pro:ctimesi}, and if $F$ in \ref{eq:chainhomotop}
is a coalgebra morphism, then the diagram \[\xymatrix{ {C\otimes I}\ar[r]^{F} \ar[d]_{\alpha_C \otimes 1} & {D}\ar[d]^{\alpha_D} \\ {P_{\mathcal{V}}(C)\otimes I}\ar[r]_-{F'}  & {P_{\mathcal{V}}D} }\]
commutes in the pointed irreducible case and the diagram \[\xymatrix{ {C\otimes I}\ar[r]^{F} \ar[d]_{\alpha_C \otimes 1} & {D}\ar[d]^{\alpha_D} \\ {L_{\mathcal{V}}(C)\otimes I}\ar[r]_-{F'}  & {L_{\mathcal{V}}D} }\]commutes
in the general case. Here $\alpha_{C}$ and $\alpha_{D}$ are classifying
maps of coalgebra structures.
\end{prop}
\begin{proof}
We will prove this in the pointed irreducible case. The general case
follows by a similar argument. The chain-homotopy between the $f_{i}$
induces\[
P_{\mathcal{V}}F:P_{\mathcal{V}}(C\otimes I)\to P_{\mathcal{V}}D\]

Now we construct the map \[
H:(P_{\mathfrak{\mathcal{V}}}C)\otimes I\to P_{\mathfrak{\mathcal{V}}}(C\otimes I)\]
 using the universal property of a cofree coalgebra and the fact that
the coalgebra structure of $(P_{\mathcal{V}}C)\otimes I$ extends
that of $P_{\mathfrak{\mathcal{V}}}C$ on both ends by condition~\ref{cond:mainoperadassumption}.
Clearly\[
P_{\mathcal{V}}F\circ H:(P_{\mathcal{V}}C)\otimes I\to P_{\mathcal{V}}D\]
is the required left-homotopy.

If we define a coalgebra structure on $C\otimes I$ using condition~\ref{cond:mainoperadassumption},
we get diagram \[\xymatrix{ {C\otimes I}\ar@{=}[r] \ar[d]_{\alpha_C \otimes 1}&{C\otimes I}\ar[r]^-{F}\ar[d]^{\alpha_{C\otimes I}} & {D}\ar[d]^{\alpha_D} \\ {P_{\mathcal{V}}(C)\otimes I}\ar[d]_{\epsilon_C \otimes 1}\ar[r]^{H} & {P_{\mathcal{V}}(C\otimes I)}\ar[r]_-{P_{\mathcal{V}}F}\ar[d]^{\epsilon_{C \otimes I}} & {P_{\mathcal{V}}D} \\  {C\otimes I}\ar@{=}[r]& {C\otimes I}& {} }\]where
$\alpha_{C\otimes I}$ is the classifying map for the coalgebra structure
on $C\otimes I$. 

We claim this diagram commutes. The fact that $F$ is a coalgebra
morphism implies that the upper right square commutes. The large square
on the left (bordered by $C\otimes I$ on all four corners) commutes
by the property of co-augmentation maps and classifying maps. The
two smaller squares on the left (i.e., the large square with the map
$H$ added to it) commute by the universal properties of cofree coalgebras
(which imply that induced maps to cofree coalgebras are uniquely determined
by their composites with co-augmentations). The diagram in the statement
of the result is just the outer upper square of this diagram, so we
have proved the claim.
\end{proof}
\begin{thm}
\label{th:freeoperadhomolinvariant}Let $\mathcal{V}$ be a cofibrant
operad whose $n^{\text{th}}$ component is $\zs{n}$-projective and
finitely generated for all $n>0$, and let \[
f:C\to D\]
be a homology equivalence of nearly free chain-complexes that are
bounded from below. Then the induced morphisms\begin{eqnarray*}
L_{\mathcal{V}}f:L_{\mathcal{V}}C & \to & L_{\mathcal{V}}D\\
M_{\mathcal{V}}f:M_{\mathcal{V}}C & \to & M_{\mathcal{V}}D\\
P_{\mathcal{V}}f:P_{\mathcal{V}}C & \to & P_{\mathcal{V}}D\\
\mathrf{F}_{\mathcal{V}}f:\mathrf{F}_{\mathcal{V}}C & \to & \mathrf{F}_{\mathcal{V}}D\end{eqnarray*}
are homology equivalences.
\end{thm}
\begin{proof}
This is a direct application of lemma~\ref{lemma:homotopyanddirectlimits},
where \[
F=H_{*}(\mathrm{suitable\,\, cofree\,\, coalgebra\,\, functor})\]

Here, we have used the fact that cofibrant operads automatically satisfy
condition~\ref{cond:mainoperadassumption}.
\end{proof}

\section{The general case\label{sec:The-general-case}}

This section states and proves theorem~\ref{th:maintheorem}.

We can relativize the definition of cofree coalgebra in definition~\ref{def:cofreecoalgebra}:

\begin{defn}
\label{def:relativecofree}Let $f:\mathcal{U}\to\mathcal{V}$ be a
morphism of operads and let $C\in\chaincat$. Any $\mathcal{V}$-coalgebra,
$A$, can be pulled back over $f$ to a $\mathcal{U}$-coalgebra,
$f^{*}A$. The \emph{relative cofree coalgebra with respect to the
morphism} $f$ \emph{and cogenerated by} $C$, denoted $L_{f}C$ solves
the universal problem:

Given any $\mathcal{V}$-coalgebra, $A$, and any morphism in $\chaincat$
$g:\forgetful{f^{*}A}\to C$, there exists a unique morphism of $\mathcal{U}$-coalgebras
$\hat{g}:f^{*}A\to L_{f}C$ that makes the diagram \[\xymatrix{ {f^*A}\ar[r]^{\hat{g}}\ar[rd]_{g} & {L_fC}\ar[d]^{\epsilon} \\ {} & {C}}\]
commute. Here, the map $\epsilon:L_{f}C\to C$ is the cogeneration
map.
\end{defn}
\begin{rem*}
These {}``not so cofree'' coalgebras are universal targets of the
subclass of $\mathcal{U}$-coalgebras that have been pulled back over
$f$. In like fashion, we can define $M_{f}C$, $P_{f}C$, and $\mathrf{F}_{f}C$.
\end{rem*}
The universal property of $L_{f}C$ immediately implies that:

\begin{prop}
\label{pro:relativecofreesub}Under the hypotheses of definition~\ref{def:relativecofree}\[
L_{f}C=\alpha_{\varepsilon}(f^{*}L_{\mathcal{V}}C)\subseteq L_{\mathcal{U}}C\]
where $\alpha_{\varepsilon}:f^{*}L_{\mathcal{V}}C\to L_{\mathcal{U}}C$
is the canonical morphism of the $\mathcal{U}$-coalgebra, $f^{*}L_{\mathcal{V}}C$
to $L_{\mathcal{U}}C$ induced by the cogeneration-projection $\varepsilon:\forgetful{f^{*}L_{\mathcal{V}}C}\to C$
(see definition~\ref{def:cofreecoalgebra}). 
\end{prop}
\begin{rem}
\label{rem:relativecofreesubrem}Corresponding statements clearly
hold for $M_{f}C$, $P_{f}C$, and $\mathrf{F}_{f}C$. The morphism
$\alpha_{\varepsilon}:f^{*}L_{\mathcal{V}}C\to L_{\mathcal{U}}C$
is not usually injective.
\end{rem}
The main idea used in theorem~\ref{th:maintheorem} is contained
in: 

\begin{lem}
\label{lem:idealannihil}Let $C\in\chaincat$ be nearly free, let
$\mainoperad$ be a projective operad that is finitely generated in
each dimension, and let $\iota:\mathcal{I}\hookrightarrow\mainoperad$
be the inclusion of an operadic ideal, inducing the map\[
\sethom(\iota,1):\sethom(\mathcal{H},C^{\otimes})\to\sethom(\mathcal{I},C^{\otimes})\]
 If $K$ is the kernel of the composite\[
\kappa:\forgetful{L_{\mainoperad}C}\xrightarrow{p}\sethom(\mathcal{H},C^{\otimes})\xrightarrow{\sethom(\iota,1)}\sethom(\mathcal{I},C^{\otimes})\]
where\[
p:C\oplus\sethom(\mathcal{H},C^{\otimes})\to\sethom(\mathcal{H},C^{\otimes})\]
is the projection, then $K$ is the pullback of a coalgebra over $\mathcal{H}/\mathcal{I}$
via the projection\[
\mathcal{H}\to\mathcal{H}/\mathcal{I}\]
 that satisfies the universal requirements for being the cofree coalgebra
$L_{\mathcal{H}/\mathcal{I}}C$.
\end{lem}
\begin{proof}
See appendix~\ref{sec:idealannihilproof} for the proof.
\end{proof}
We can prove corresponding statements for the truncated and pointed-irreducible
cofree coalgebras:

\begin{cor}
\label{cor:annihextensions}Under the hypotheses of lemma~\ref{lem:idealannihil},
if $M$ is the kernel of the composite\[
\forgetful{M_{\mainoperad}C}\xrightarrow{p}\sethom(\mathcal{H},C^{\otimes})\xrightarrow{\sethom(\iota,1)}\sethom(\mathcal{I},C^{\otimes})\]
where $k=0$ if $\mainoperad$ is unital and $1$ otherwise, then
$M=\im M_{\mainoperad/\mathcal{I}}C$ in $M_{\mainoperad}C$ under
the natural map induced by the projection $\mainoperad\to\mainoperad/\mathcal{I}$.

If $\mainoperad$is a unital operad and $P$ is the kernel of the
composite

\[
\forgetful{P_{\mainoperad}C}\xrightarrow{p}\sethom(\mathcal{H},C^{\otimes})\xrightarrow{\sethom(\iota,1)}\sethom(\mathcal{I},C^{\otimes})\]
then $P=\im P_{\mainoperad/\mathcal{I}}C\subseteq P_{\mainoperad}C$.
If $\mathrf{F}$ is the kernel of the composite\[
\forgetful{\mathrf{F}_{\mainoperad}C}\xrightarrow{p}\sethom(\mathcal{H},C^{\otimes})\xrightarrow{\sethom(\iota,1)}\sethom(\mathcal{I},C^{\otimes})\]
and $\mainoperad/\mathcal{I}$ is a unital operad, then $\mathrf{F}=\im\mathrf{F}_{\mainoperad/\mathcal{I}}C\subseteq\mathrf{F}_{\mainoperad}C$.
\end{cor}
\begin{proof}
The proof of lemma~\ref{lem:idealannihil} does not use any specific
property of $L_{\mainoperad}C$ other than the facts that
\begin{enumerate}
\item it is a coalgebra that is a submodule of $\sethom(\mathcal{H},C^{\otimes})$
\item its coproduct is dual to the compositions of $\mainoperad$
\item it is cofree in a suitable context
\end{enumerate}
It is only necessary to remark that the fact that $\mainoperad/\mathcal{I}$
is unital implies that $\eta(1)\notin\mathcal{I}_{1}$ so that the
basepoint of $P_{\mainoperad}C$ and $\mathrf{F}_{\mainoperad}C$
lie in $P$ and $\mathrf{F}$, respectively.
\end{proof}
Now we define functoriality of cofree coalgebras with respect to operad-morphisms:

\begin{thm}
\label{th:maintheorem}Let $f:\mathcal{I}\hookrightarrow\mathcal{H}$
be the inclusion of an operadic ideal with $\mathcal{H}$ a projective
free operad, $\mathcal{V}=\mathcal{H}/\mathcal{I}$ a projective operad,
and with canonical projection $p:\mathcal{H}\to\mathcal{H}/\mathcal{I}=\mathcal{V}$.
In addition, let $C\in\chaincat$ be nearly free. Then the kernels
of \[
\left\{ \begin{array}{c}
\alpha_{\varepsilon}:f^{*}L_{\mathcal{H}}C\to L_{\mathcal{I}}C\\
\alpha_{\varepsilon}:f^{*}M_{\mathcal{H}}C\to M_{\mathcal{I}}C\\
\alpha_{\varepsilon}:f^{*}P_{\mathcal{H}}C\to P_{\mathcal{I}}C\\
\alpha_{\varepsilon}:f^{*}\mathrf{F}_{\mathcal{H}}C\to\mathrf{F}_{\mathcal{I}}C\end{array}\right\} \]
(see definition~\ref{def:cofreecoalgebra} and proposition~\ref{pro:relativecofreesub}
for an explanation of the notation $\alpha_{\varepsilon}$) are\[
\left\{ \begin{array}{c}
\forgetful{p^{*}L_{\mathcal{V}}C}/C\\
\forgetful{p^{*}M_{\mathcal{V}}C}/C\\
\forgetful{p^{*}P_{\mathcal{V}}C}/\ring\\
\forgetful{p^{*}\mathrf{F}_{\mathcal{V}}C}/\ring\end{array}\right\} \]
 respectively. If\[
W_{*}C=\left\{ \begin{array}{c}
L_{*}C\\
M_{*}C\\
P_{*}C\\
\mathrf{F}_{*}C\end{array}\right\} \]
then $W_{f}C\subseteq W_{\mathcal{I}}C$ has the structure of an $\mathcal{H}$-coalgebra.
This coalgebra structure induces an $\mathcal{H}$-coalgebra morphism
$\hat{f}:W_{f}C\to W_{\mathcal{H}}C$ that is a right inverse to $f^{*}$.
This, in turn, induces a splitting of underlying chain-complexes\begin{equation}
\forgetful{W_{\mathcal{H}}C}\cong\forgetful{W_{\mathcal{V}}C}/C\oplus\forgetful{W_{f}C}\label{eq:splittingchain}\end{equation}
If $\mainoperad$is finitely generated in each dimension any homology
equivalence\[
C\to C'\]
of nearly free modules that are bounded from below induces a homology
equivalence\[
W_{\mathcal{V}}C\to W_{\mathcal{V}}C'\]

\end{thm}
\begin{rem*}
Note that $W_{f}C=\alpha_{\varepsilon}(f^{*}W_{\mathcal{H}}C)\subseteq L_{\mathcal{I}}C$,
by proposition~\ref{pro:relativecofreesub} and remark~\ref{rem:relativecofreesubrem}.

This result's key ideas can be summarized as follows:
\begin{enumerate}
\item An operad morphism $f:\mathcal{U}\to\mathcal{V}$ induces a $\mathcal{U}$-coalgebra
morphism\[
f^{*}:W_{\mathcal{V}}C\to W_{\mathcal{U}}C\]
whose kernel is a priori a \emph{coideal} (see definition~\ref{def:coideal}).
\item In the special case where $f:\mathcal{I}\to\mainoperad$ is the inclusion
of an operadic ideal, the kernel of the induced map\[
f^{*}:W_{\mathcal{H}}C\to W_{\mathcal{I}}C\]
is a full-fledged \emph{$\mainoperad$-coalgebra} (this is the main
thrust of appendix~\ref{sec:idealannihilproof}), hence the image
of $f^{*}$ is also an $\mainoperad$-coalgebra (being the quotient
of two such).
\item The universal property of a cofree coalgebra implies the existence
of a \emph{unique} map\[
c:f^{*}(W_{\mainoperad}C)\to W_{\mainoperad}C\]
\emph{splitting} $f^{*}$.
\item This universal property also implies that the kernel is the pullback
of $W_{\mainoperad/\mathcal{I}}C$.
\item The homology invariance of $W_{\mainoperad}C$ --- established in
theorem~\ref{th:freeoperadhomolinvariant} --- implies that of $W_{\mainoperad/\mathcal{I}}C$.
\end{enumerate}
\end{rem*}
\begin{proof}
We will prove this in the case where $W_{*}C=L_{*}C$. The other cases
follow by similar arguments.

Lemma~\ref{lem:idealannihil} implies that the kernel, $Z$, of\[
\alpha_{\varepsilon}:f^{*}L_{\mathcal{H}}C\to L_{\mathcal{I}}C\]
 is the pullback of a coalgebra over $\mathcal{V}=\mathcal{H}/\mathcal{I}$.
Here, the fact that $\mathcal{I}$ is an operadic ideal implies that
$Z$ is a \emph{sub-coalgebra} rather than a mere coideal --- indeed,
it is $p^{*}L_{\mathcal{V}}C$. The subcoalgebra $Z\oplus C\subseteq L_{\mathcal{H}}C$
(where $C$ is equipped with a coproduct that is identically $0$)
is also the pullback of a coalgebra over $\mathcal{V}$ and has the
universal property of $p^{*}L_{\mathcal{V}}C$ so $Z\oplus C=p^{*}L_{\mathcal{V}}C$.
Consider\[
\epsilon_{\mathcal{V}}:p^{*}L_{\mathcal{V}}C\to C\]
where $C$ is regarded as a $\mathcal{V}$-coalgebra whose coproduct
is \emph{identically zero}. The kernel of $\epsilon_{\mathcal{V}}$
will be a \emph{coideal} in $p^{*}L_{\mathcal{V}}C$ (see definition~\ref{def:coideal})
whose underlying chain complex is isomorphic to $\forgetful{p^{*}L_{\mathcal{V}}C}/C$
(since $C\subset p^{*}L_{\mathcal{V}}C$ is a direct summand as a
chain complex and as a coalgebra).

We claim that $\ker\epsilon_{\mathcal{V}}$ is \emph{also} a coideal
in $L_{\mathcal{H}}C$. Consider the diagram \[\xymatrix{{\ker \epsilon_{\mathcal{V}}}\ar[d]\ar[dr] &  {}\\ {\sethom(\mathcal{H},(p^*L_{\mathcal{V}}C)^\otimes)}\ar[d]\ar[r] & {\sethom(\mathcal{H},(L_{\mathcal{H}}C)^\otimes)}\ar[d] \\{\sethom(\mathcal{H},(p^*L_{\mathcal{V}}C/\ker \epsilon_{\mathcal{V}})^\otimes)}\ar[r] & {\sethom(\mathcal{H},(L_{\mathcal{H}}C/\ker \epsilon_{\mathcal{V}})^\otimes)}}\]where
the maps from $\ker\epsilon_{\mathcal{V}}$ are the structure maps
of $p^{*}L_{\mathcal{V}}C$ and $L_{\mathcal{H}}C$ and the remaining
downward maps are induced by projection to the quotient. The upper
triangle commutes since $p^{*}L_{\mathcal{V}}C=Z\oplus C$ is a sub-coalgebra
of $L_{\mathcal{H}}C$. The remaining squares commute by naturality
of projection to the quotient. 

The composite of the vertical maps on the left is $0$ because $\ker\epsilon_{\mathcal{V}}$
is a coideal in $p^{*}L_{\mathcal{V}}C$ --- see definition~\ref{def:coideal}.
The commutativity of the diagram implies that the composite of the
vertical maps on the right is also $0$, so $\ker\epsilon_{\mathcal{V}}$
is a coideal in $L_{\mathcal{H}}C$.

It follows that the quotient \[
L_{f}C=L_{\mathcal{H}}C/\ker\epsilon_{\mathcal{V}}\subseteq L_{\mathcal{I}}C\]
 is an $\mathcal{H}$-coalgebra equipped with a canonical cogeneration
(chain-)map\[
\epsilon_{f}:\forgetful{L_{f}C}\to C\]

This chain-map and the universal property of the cofree coalgebra
$L_{\mathcal{H}}C$ implies the existence of a coalgebra morphism\[
\hat{f}:L_{f}C\to L_{\mathcal{H}}C\]
The composite of this with $\alpha_{\varepsilon}$ is a morphism that
covers the identity map of $C$ --- which must \emph{be} the identity
map of $L_{f}C\subseteq L_{\mathcal{I}}C$ due to the uniqueness of
induced maps to cofree coalgebras. Consequently, $\hat{f}$ splits
$\alpha_{\varepsilon}$ and induces the splitting of chain-complexes
in equation~\ref{eq:splittingchain}.

The final statements follows from lemma~\ref{lem:idealannihil} and
the fact that every operad is the surjective image of some free operad.
So the splitting in equation~\ref{eq:splittingchain} exists for
any $\mathcal{V}$ and suitable free operad. This splitting induces
a corresponding splitting in homology\[
H_{*}(\forgetful{L_{\mathcal{H}}C})\cong H_{*}(\forgetful{L_{\mathcal{V}}C}/C)\oplus H_{*}(\forgetful{L_{f}C})\]
The statement about homology invariance of $L_{\mathcal{V}}C$ follows
from theorem~\ref{th:freeoperadhomolinvariant} and the fact that
a direct summand of an isomorphism is an isomorphism.
\end{proof}

\begin{cor}
\label{cor:mainresult} Let $\ring$ be a field or $\integers$ and
let $\mathcal{V}=\{\mathcal{V}(n)\}$ be an operad such that $\mathcal{V}(n)$
is $\ring S_{n}$-projective and finitely generated in each dimension
for all $n>0$. If{\[
W_{\mathcal{V}}C=\left\{ \begin{array}{c}
L_{\mathcal{V}}C\\
M_{\mathcal{V}}C\\
P_{\mathcal{V}}C\\
\mathrf{F}_{\mathcal{V}}C\end{array}\right\} \]
} and\[
f:C\to D\]
is a homology equivalence of nearly free chain complexes (see definition~\ref{def:nearlyfree})
that are bounded from below, then the induced map\[
W_{\mathcal{V}}f:W_{\mathcal{V}}C\to W_{\mathcal{V}}D\]
is a homology equivalence. 
\end{cor}
\begin{proof}
Given $\mathcal{V}$ satisfying the hypotheses, let $\mainoperad$
be the free operad generated by the components of $\mathcal{V}$.
It will satisfy the hypotheses of theorem~\ref{th:maintheorem} and
there will exist a canonical surjection of operads\[
\mainoperad\to\mathcal{V}\]
whose kernel is an operadic ideal.
\end{proof}
\appendix

\section{Nearly free modules\label{sec:nearlyfree}}

In this section, we will explore the class of nearly free $\integers$-modules
--- see definition~\ref{def:nearlyfree}. We show that this is closed
under the operations of taking direct sums, tensor products, countable
products and cofree coalgebras. It appears to be fairly large, then,
and it would be interesting to have a direct algebraic characterization.

Clearly a module must be torsion-free (hence flat) to be nearly free.
The converse is not true, however: $\mathbb{Q}$ is flat but \emph{not}
nearly free.

The definition immediately implies that:

\def\ring{\integers}

\begin{prop}
\label{prop:nfsubmodule}Any submodule of a nearly free module is
nearly free.
\end{prop}
Nearly free modules are closed under operations that preserve free
modules: 

\begin{prop}
\label{prop:nfsumotimes}Let $M$ and $N$ be $\integers$-modules.
If they are nearly free, then so are $M\oplus N$ and $M\otimes N$.

Infinite direct sums of nearly free modules are nearly free.
\end{prop}
\begin{proof}
If $F\subseteq M\oplus N$ is countable, so are its projections to
$M$ and $N$, which are free by hypothesis. It follows that $F$
is a countable submodule of a free module.

The case where $F\subseteq M\otimes N$ follows by a similar argument:
The elements of $F$ are finite linear combinations of monomials $\{ m_{\alpha}\otimes n_{\alpha}\}$
--- the set of which is countable. Let\begin{eqnarray*}
A & \subseteq & M\\
B & \subseteq & N\end{eqnarray*}
 be the submodules generated, respectively, by the $\{ m_{\alpha}\}$
and $\{ n_{\alpha}\}$. These will be countable modules, hence $\integers$-free.
It follows that \[
F\subseteq A\otimes B\]
is a free module.

Similar reasoning proves the last statement, using the fact that any
direct sum of free modules is free.
\end{proof}
\begin{prop}
\label{prop:nfprodfree}Let $\{ F_{n}\}$ be a countable collection
of $\integers$-free modules. Then\[
\prod_{n=1}^{\infty}F_{n}\]
is nearly free.
\end{prop}
\begin{proof}
In the case where $F_{n}=\integers$ for all $n$\[
B=\prod_{n=1}^{\infty}\integers\]
is the Baer-Specker group, which is well-known to be nearly free ---
see \cite{Baer:1937}, \cite[vol. 1, p. 94 Theorem 19.2]{Fuchs:1970},
and\cite{Blass-Gobel:1996}. It is also well-known \emph{not} to be
$\integers$-free --- see \cite{Baer-Specker-nonfree} or the survey
\cite{Baer-Specker-survey}.

In the general case, \[
\prod_{n=1}^{\infty}F_{n}\]
is a direct sum of copies of $B$, which is nearly free by proposition~\ref{prop:nfsumotimes}.
\end{proof}
\begin{cor}
\label{cor:nfprodnf}Let $\{ N_{k}\}$ be a countable set of nearly
free modules. Then\[
\prod_{k=1}^{\infty}N_{k}\]
is also nearly free.
\end{cor}
\begin{proof}
Let\[
F\subset\prod_{k=1}^{\infty}N_{k}\]
 be countable. If $F_{k}$ is its projection to factor $N_{k}$, then
$F_{k}$will be countable, hence free. It follows that\[
F\subset\prod_{k=1}^{\infty}F_{k}\]
 and the conclusion follows from proposition~\ref{prop:nfprodfree}. 
\end{proof}
\begin{cor}
\label{cor:nfhom1}Let $A$ be nearly free and let $F$ be $\integers$-free
of countable rank. Then\[
\homz(F,A)\]
is nearly free.
\end{cor}
\begin{proof}
This follows from corollary~\ref{cor:nfprodnf} and the fact that
\[
\homz(F,A)\cong\prod_{k=1}^{\mathrm{rank}(F)}A\]

\end{proof}
\begin{cor}
\label{cor:nfprodfn}Let $\{ F_{n}\}$ be a sequence of $\zs{n}$-projective
modules and and let $A$ be nearly free. Then\[
\prod_{n=1}^{\infty}\homzs{n}(F_{n},A^{n})\]
 is nearly free.
\end{cor}
\begin{proof}
This is a direct application of the results of this section and the
fact that \[
\homzs{n}(F_{n},A^{n})\subseteq\homz(F_{n},A^{n})\subseteq\homz(\hat{F}_{n},A^{n})\]
where $\hat{F}_{n}$ is a $\zs{n}$-free module of which $F_{n}$
is a direct summand.
\end{proof}
\begin{thm}
\label{thm:nfcofreenf}Let $C$ be a nearly free $\integers$-module
and let $\mathcal{V}$ be an operad whose $n^{\text{th}}$ component
is $\zs{n}$-projective and finitely generated for all $n$. Then\begin{eqnarray*}
 & \forgetful{L_{\mathcal{V}}C}\\
 & \forgetful{M_{\mathcal{V}}C}\\
 & \forgetful{P_{\mathcal{V}}C}\\
 & \forgetful{\mathrf{F}_{\mathcal{V}}C}\end{eqnarray*}
are all nearly free.
\end{thm}
\begin{proof}
This follows from theorem~\ref{th:cofreedirectlimit} which states
that all of these are submodules of \[
\prod_{n=1}^{\infty}\homzs{n}(\mathcal{V}_{n},A^{n})\]
and the fact that near-freeness is inherited by submodules.
\end{proof}

\section{The relative coendomorphism operad of the unit interval\label{sec:relativecoendomorphism}}

Our main result is:

\begin{prop}
\label{pro:computerelativecoendunit}If $I$ is the unit interval
(see definition~\ref{def:unitinterval}), its relative coendomorphism
operad (see definition~\ref{def:coend}) is given by\[
\coend(I;\{\integers\cdot p_{0},\integers\cdot p_{1}\})=\mathfrak{S}_{0}\]
defined in Example~\ref{example:frakS0}.
\end{prop}
\begin{proof}
We must compute homomorphisms\[
g:I\to I^{n}\]
that send the endpoints $\{ p_{0},p_{1}\}$ to the subcomplex of $I^{n}$
generated by tensor products of the endpoints --- i.e. \[
\integers p_{0}\otimes\cdots\otimes p_{0}\oplus\integers p_{1}\otimes\cdots\otimes p_{1}\]
Both of these subcomplexes (of $I$ and $I^{n}$) are concentrated
in dimension $0$, which implies that all of our maps must be of degree
zero.

It follows that all components of $\coend(I;\{\ring\cdot p_{0},\ring\cdot p_{1}\})$
are concentrated in dimension $0$. Chain-maps of $I$ are determined
by where they send the 1-dimensional element, $q$. Thus we want chain-maps\[
g:I\to I^{n}\]
 with $\partial(g(q))=p_{1}\otimes\cdots\otimes p_{1}-p_{0}\otimes\cdots\otimes p_{0}$
($n$ factors in each term).

We use a {}``geometric argument.'' Consider the unit cube in $\mathbb{R}^{n}$
with coordinates \[
0\le x_{i}\le1\]
for $i=1,\dots,n$. Regard the edges of this as 1-simplices and the
vertices as 0-simplices. Chains with the required property correspond
to sequences of these 1-simplices forming paths along the edges of
the cube from $(0,\dots0)$ to $(1,\dots,1)$.

We claim there are exactly $n!$ such paths and they are linearly
independent chains in $C(I^{n})_{1}$. To construct a path, one must
travel 1 unit in the $x_{i}$ direction, then 1 unit in the $x_{i'}$
direction, with $i'\ne i$, and so on. One represents this by a list
of $n$ distinct integers between 1 and $n$:\[
(i,i',\dots)\]

Such lists clearly correspond to permutations $\sigma\in S_{n}$:\[
(\sigma(1),\sigma(2),\dots,\sigma(n))\]
 Let $\{ v_{0},\dots,v_{n}\}$ be coordinates of the vertices one
encounters during this process with $v_{0}=(0,\dots,0)$ and $v_{n}=(1,\dots,1)$. 

Since $v_{k+1}-v_{k}$ determines the direction one went in the $k^{\text{th}}$step
(and since each path travels in a direction taken by no other in \emph{some}
step), it follows that each path has a \emph{vertex} not contained
in any other. This implies that each path also has a \emph{1-simplex}
not contained in any other. Consequently the paths represent linearly
independent chains of $C(I^{n})_{1}$.

It is also clear that the symmetric group permutes these $n!$ paths
by permuting coordinate axes. This demonstrates a natural equality\[
\coend(I;\{\integers\cdot p_{0},\integers\cdot p_{1}\})([n])=\integers S_{n}\]

\end{proof}

\section{Homotopy and direct limits\label{sec:homotopyanddirect}}

This section's main result may be summed up by the phrase

\begin{quotation}
{}``A homotopy functor that commutes with direct limits is a homology
functor of nearly free complexes.''
\end{quotation}
Let $\homotopycat(\integers)$ denote the \emph{chain-homotopy category}
of $\integers$-chain-complexes --- compare to the notation in \S~20.4
of \cite{Weibel:homological-algebra}. Objects in this category are
chain-homotopy equivalence classes of chain complexes (not necessarily
torsion free) and chain-homotopic morphisms are equivalent. 

We have the related category, $\mathbf{D}(\integers)$ --- essentially
the Verdier derived category of $\integers$. Its objects are chain-complexes
where homology equivalent complexes are considered equivalent (the
Verdier derived category considered cochain complexes).

We also consider the subcategory $\homotopcellular\subseteq\homotopycat(\integers)$
of cellular chain-complexes --- Exercise~10.4.5 of \cite{Weibel:homological-algebra}.
These are chain complexes \[
C=\bigcup_{i=1}^{\infty}C_{i}\]
where $C_{i+1}/C_{i}$ is $\integers$-free and has vanishing differential.
Clearly, any $\integers$-free chain complex that is bounded from
below is in this category.

We use the following well-known properties of $\homotopcellular$
and $\homotopycat(\integers)$:

\begin{enumerate}
\item If $C\in\homotopcellular$ and $A\in\chaincat$ is acyclic, then every
map\[
C\to A\]
is nullhomotopic.
\item If $C\in\homotopcellular$ and\[
f:A\to B\]
is a homology equivalence in $\chaincat$, then \[
f_{*}:\homotophom(C,A)\xrightarrow{\cong}\homotophom(C,B)\]
is an isomorphism.
\item If $C,D\in\homotopcellular$ and\[
f:C\to D\]
is a homology equivalence then $f$ is also a homotopy equivalence.
\end{enumerate}
Our main result is:

\begin{lem}
\label{lemma:homotopyanddirectlimits}Let $F:\chaincat\to\mathrm{mod}-\integers$
be a functor such that
\begin{enumerate}
\item whenever $\{ C_{\alpha}\}$ is a direct system of cellular complexes
in $\chaincat$, \[
F(\dlimit C_{\alpha})=\dlimit F(C_{\alpha})\]

\item $F$ factors through the natural quotient $\chaincat\to\homotopycat_{\mathrm{cell}}$
(i.e., $F$ is a homotopy functor). 
\end{enumerate}
If\[
f:C\to D\]
is a homology equivalence of nearly free chain-complexes that are
bounded from below, then\[
F(f):F(C)\to F(D)\]
is an isomorphism.

\end{lem}
\begin{rem*}
This essentially says
\end{rem*}
\begin{quotation}
A homotopy functor that commutes with direct limits is a homology
functor of nearly free complexes.
\end{quotation}
\begin{proof}
The conclusion is already known to be true if $C$ and $D$ are in
$\homotopcellular$ because then they are homotopy equivalent and
$F$ is assumed to be a homotopy functor. 

In the general case, let\begin{eqnarray*}
C & = & \dlimit C_{\alpha}\\
D & = & \dlimit D_{\alpha}\end{eqnarray*}
 where $C_{\alpha}$ and $D_{\alpha}$ are countable chain-complexes.
This is possible by proposition~\ref{pro:directlimitfintelygenerated}. 

In addition, assume $D_{\alpha}=f(C_{\alpha})$ --- since homomorphic
images of countable complexes are countable. There may be other $D_{\alpha'}$
not in the image of any of the $C_{\alpha}$. Our hypotheses imply
that\begin{eqnarray*}
F(C) & = & \dlimit F(C_{\alpha})\\
F(C) & = & \dlimit F(D_{\alpha})\end{eqnarray*}

Let \begin{eqnarray*}
c_{\alpha}:C_{\alpha} & \to & C\\
c_{\alpha,\beta}:C_{\alpha} & \to & C_{\beta}\\
d_{\alpha}:D_{\alpha} & \to & D\\
d_{\alpha,\beta}:D_{\alpha} & \to & D_{\beta}\end{eqnarray*}
be the inclusions.

The properties of $\homotopcellular$ imply the commutativity of

\begin{equation}\label{dia:commute1}\xymatrix{{\homotophom(D_\alpha,C)}\ar[r]^{\cong}\ar[d]_{\homotophom(f_\alpha,1)} & {\homotophom(D_\alpha,D)}\ar[d]^{\homotophom(f_\alpha,1)} \\ {\homotophom(C_\alpha,C)}\ar[r]^{\cong} & {\homotophom(C_\alpha,D)}}\end{equation}
in the case where $D_{\alpha}=f(C_{\alpha})$, and\begin{equation}\label{dia:commute2}\xymatrix{{\homotophom(D_\beta,C)}\ar[r]^{\cong}\ar[d]_{\homotophom(d_{\alpha,\beta},1)} & {\homotophom(D_\beta,D)}\ar[d]^{\homotophom(d_{\alpha,\beta},1)} \\ {\homotophom(D_\alpha,C)}\ar[r]^{\cong} & {\homotophom(D_\alpha,D)}}\end{equation}whenever
$D_{\alpha}\subseteq D_{\beta}$. 

Let $h_{\alpha}\in\homotophom(D_{\alpha},C)$ map to $d_{\alpha}\in\homotophom(D_{\alpha},D)$
under the isomorphism above. They are maps\[
h_{\alpha}:D_{\alpha}\to C\]
that are well-defined up to \emph{homotopy}.

Diagrams \ref{dia:commute1} and \ref{dia:commute2} implies the \emph{homotopy}
commutativity of the diagrams

$$
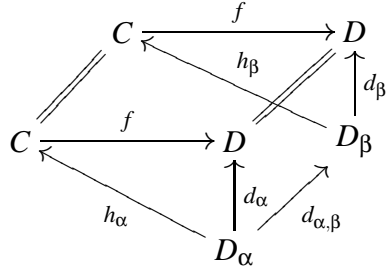
\xymatrix{{}&{C}\ar[rr]^f&{}&{D}\\{C}\ar[rr]^f\ar@{=}[ur]&{}&{D}\ar@{=}[ru]&{D_\beta}\ar[ull]_{h_\beta}\ar[u]_{d_\beta}\\{}&{}&{D_\alpha}\ar[ull]^{h_\alpha}\ar[u]_{d_\alpha}\ar[ru]_{d_{\alpha,\beta}}&{}}$$and

 
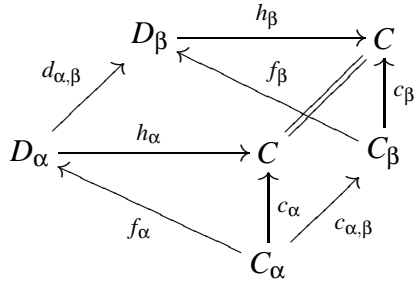
$$\xymatrix{{}&{D_\beta}\ar[rr]^{h_\beta}&{}&{C}\\{D_\alpha}\ar[rr]^{h_\alpha}\ar[ur]^{d_{\alpha,\beta}}&{}&{C}\ar@{=}[ru]&{C_\beta}\ar[ull]_{f_\beta}\ar[u]_{c_\beta}\\{}&{}&{C_\alpha}\ar[ull]^{f_\alpha}\ar[u]_{c_\alpha}\ar[ru]_{c_{\alpha,\beta}}&{}}$$whenever
$D_{\alpha}\subseteq D_{\beta}$. 

The fact that $F$ is a homotopy functor implies the \emph{exact}
commutativity of the diagrams 
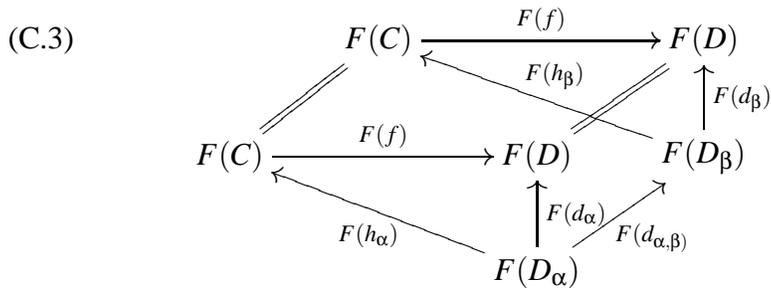
\begin{equation}\label{dia:commute3}\xymatrix{{}&{F(C)}\ar[rr]^{F(f)}&{}&{F(D)}\\{F(C)}\ar[rr]^{F(f)}\ar@{=}[ur]&{}&{F(D)}\ar@{=}[ru]&{F(D_\beta)}\ar[ull]_{F(h_\beta)}\ar[u]_{F(d_\beta)}\\{}&{}&{F(D_\alpha)}\ar[ull]^{F(h_\alpha)}\ar[u]_{F(d_\alpha)}\ar[ru]_{F(d_{\alpha,\beta)}}&{}}\end{equation}and

 \begin{equation}\label{dia:commute4}\xymatrix{{}&{F(D_\beta)}\ar[rr]^{F(h_\beta)}&{}&{F(C)}\\{F(D_\alpha)}\ar[rr]^{F(h_\alpha)}\ar[ur]^{F(d_{\alpha,\beta})}&{}&{F(C)}\ar@{=}[ru]&{F(C_\beta)}\ar[ull]_{F(f_\beta)}\ar[u]_{F(c_\beta)}\\{}&{}&{F(C_\alpha)}\ar[ull]^{F(f_\alpha)}\ar[u]_{F(c_\alpha)}\ar[ru]_{F(c_{\alpha,\beta})}&{}}\end{equation}when
they are well-defined.

Diagrams \ref{dia:commute3} (for all values of $\alpha$) imply the
existence of a map\[
h:\dlimit F(D_{\alpha})=F(D)\to F(C)\]
that is a right-inverse to $F(f):F(C)\to F(D)$, and diagrams \ref{dia:commute4}
imply that it is also a left-inverse.
\end{proof}

\section{\label{sec:idealannihilproof}Proof of Lemma~\ref{lem:idealannihil}}

Compare the following definition with definition~3.1 in \cite{Smith-cofree}:

\begin{defn}
\label{def:addtree}Let $k$ be $0$ or $1$. Define $\mathrf{P}(k)$
to be the set of finite sequences $\{ u_{1},\dots,u_{m}\}$ of elements
each of which is either a $\bullet$-symbol or an integer $\ge k$.

Given a sequence $\mathbf{u}\in\mathrf{P}(k)$, let $\slength{\mathbf{u}}$
denote the length of the sequence.
\end{defn}
\begin{rem}
Throughout the rest of this section, we set $k=0$ if $\mainoperad$
is unital and $k=1$ otherwise.

If $\mathrf{P}_{k}(n)$ is as defined in {definition~3.1
of \cite{Smith-cofree}}, it is not hard to see that \[
\mathrf{P}(k)=\bigcup_{n=1}^{\infty}\mathrf{P}_{k}(n)\]
 
\end{rem}

\begin{defn}
\label{def:generalizedcomps}Let $\mathcal{V}$ be an operad and let
$\mathbf{u}=\{ u_{1},\dots,u_{m}\}\in\mathrf{P}(k)$, where we impose
no condition on $k$. We define the \emph{generalized composition}
with respect to $\mathbf{u}$, denoted $\gamma_{\mathbf{u}}$, by
\begin{multline*}
\gamma(\mathbf{u})=\circ_{u_{m}}(\circ_{u_{m-1}}\otimes1)\cdots(\circ_{u_{1}}\otimes\cdots\otimes1)\circ\bigotimes_{\mathbf{u}}\iota_{j}\\
:\mathcal{V}(\mathbf{u})\otimes\bigotimes_{\mathbf{u}}\mathcal{V}(u_{j})\to\mathcal{V}(\bigsqcup_{i=1}^{m}u_{i})=\mathcal{V}(\mathfrak{g}(\mathbf{u}))\end{multline*}
where we follow the convention that
\begin{enumerate}
\item $\mathcal{V}(\{\bullet\})=\ring$,
\item $\circ_{\{\bullet\}}=\circ_{\{ x\}}\circ(\eta_{\{ x\}}\otimes1):\ring\otimes\mathcal{V}(s)=\mathcal{V}(s)\to\mathcal{V}(s\setminus\{\bullet\}\sqcup\{ x\})$,
where $\{ x\}$ is a singleton set \emph{not} containing the distinguished
element $\bullet$.
\end{enumerate}
\end{defn}
\begin{rem*}
See definition~\ref{def:setfmodules} for the definition of $\mathfrak{g}(\mathbf{u})$.

If $\mathbf{v}=\{ u_{k_{1}},\dots,u_{k_{t}}\}\subset\{ u_{1},\dots,u_{m}\}$
is the subset of non-$\bullet$ sets, then $\gamma_{\mathbf{u}}$
is a map\[
\gamma(\mathbf{u}):\mathcal{V}(\mathbf{u})\otimes\bigotimes_{\mathbf{v}}\mathcal{V}(u_{j})\to\mathcal{V}(\bigsqcup_{i=1}^{m}u_{i})=\mathcal{V}(\mathfrak{g}(\mathbf{u}))\]

If $u\in\mathrf{P}(k)$ with $x\in u$, then $u\sqcup_{x}x$ represents
$(u\setminus x)\sqcup x$ --- we have removed $x$ from $u$ and then
added the \emph{contents of} $x$ to $u$. For this notation to make
any sense, $x$ must be a set, not an atomic element. Definition~\ref{def:generalizedcomps}
to make any sense, the elements of $\mathbf{u}$ must all be sets
and the result of carrying out this operation on all of the elements
of $\mathbf{u}$ will be the {}``flattened form'' of $\mathbf{u}$
or $\mathfrak{g}(\mathbf{u})$.
\end{rem*}
Recall that $\mainoperad$ is a projective operad with operadic ideal
$\iota:\mathcal{I}\hookrightarrow\mainoperad$ and $K$ is the kernel
of the composite 

\[
\kappa:\forgetful{L_{\mainoperad}C}\xrightarrow{p}\sethom(\mathcal{H},C^{\otimes})\xrightarrow{\sethom(\iota,1)}\sethom(\mathcal{I},C^{\otimes})\]
We will show that the coalgebra structure of $L_{\mainoperad}C$ induces
a coalgebra structure on $K$ that makes it a coalgebra over $\mainoperad/\mathcal{I}$
--- pulled back over the projection\[
\mainoperad\to\mainoperad/\mathcal{I}\]
It will then turn out to inherit the {}``cofreeness'' of $L_{\mainoperad}C$
as well.

\begin{prop}
\label{pro:kerunorderedspan}Let $X\in\setf$, $x\in X$ and $\{ f_{y}:V_{y}\to U_{y}\}$
be as in definition~\ref{def:unorderedspans}. Then \[
\bigcap_{x\in X}\ker\bigotimes_{X,x}(1,f_{x})=\bigotimes_{X}\ker f_{x}\]

\end{prop}
\begin{proof}
The flatness of all the underlying modules implies that \[
\ker\bigotimes_{X,x}(1,f_{x})=\bigotimes_{X,x}(V,\ker f_{x})=V_{x_{1}}\otimes\cdots\otimes V_{x_{k}}\otimes\ker f_{x}\otimes V_{x_{k+2}}\otimes\cdots\otimes V_{x_{t}}\]
and the conclusion follows.
\end{proof}
Clearly, $K$ inherits a map\[
a:K\to\sethom(\mathcal{H},C^{\otimes})\]
from its inclusion into $L_{\mainoperad}C$. We must show that its
image actually lies in\[
\sethom(\mathcal{H},K^{\otimes})\subseteq\sethom(\mathcal{H},C^{\otimes})\]
We make use of the fact that the structure-map of $L_{\mainoperad}C$
is dual to the compositions of the operad $\mainoperad$ and that
$\mathcal{I}$ is an operadic ideal. 

The construction of $L_{\mainoperad}C$ in \cite{Smith-cofree} implies
that the diagram

\begin{equation} \xymatrix@C+30pt{{} & {\displaystyle\sethom(\mainoperad,(L_{\mainoperad}C)^\otimes)} \ar@{^{(}->}[d]^-{y}\\ {L_{\mainoperad}C\vphantom{\displaystyle\prod_{\mathbf{u}\in\mathrf{P}_k(n) }}}\ar@<1ex>[r]_-{g}\ar@<.4ex>[ru]^{\alpha} & {\displaystyle\!\!\!\!\!\prod_{ \mathbf{u}\in\mathrf{P}(k) }\!\!\!\!\!\homz(\mainoperad({\mathbf{u}})\otimes \bigotimes_{u_i\in\mathbf{u}}\mainoperad(u_{i}),C^{\mathfrak{g}(\mathbf{u})})}} \label{eq:bigdiag}\end{equation}commutes.
This is just diagram~3.2 in \cite{Smith-cofree}, where:

\begin{enumerate}
\item $\alpha$ is the adjoint structure map.
\item $\kappa:L_{\mainoperad}C\hookrightarrow\sethom(\mainoperad,C^{\otimes})$
is the inclusion (see theorem~\ref{th:cofreedirectlimit}).
\item $g=\left(\prod_{\mathbf{u}\in\mathrf{P}(k)}c(\mathbf{u})\right)\circ\kappa$
and the $c(\mathbf{u})$ are defined by \begin{multline*}
c(\mathbf{u})=\homz(\gamma(\mathbf{u}),1):\sethomn{\mathbf{u}}(\mathcal{V}(\mathfrak{g}(\mathbf{u})),C^{\mathfrak{g}(\mathbf{u})})\\
\to\homz(\mathcal{V}(\mathbf{u})\otimes\bigotimes_{\mathbf{u}}\mathcal{V}(u_{i}),C^{\mathfrak{g}(\mathbf{u})})\end{multline*}
 --- the dual of the generalized structure-map \[
\gamma(\mathbf{u}):\mathcal{V}(\mathbf{u})\otimes\bigotimes_{\mathbf{u}}\mathcal{V}(u_{i})\to\mathcal{V}(\mathfrak{g}(\mathbf{u}))\]
 from definition~\ref{def:generalizedcomps}. We assume that $\mathcal{V}(\bullet)=\ring$
and $C^{\bullet}=C$ so that $\sethom(\mathcal{V}(\bullet),C^{\bullet})=C$.
\item \label{list:yudef}if $P=\sethom(\mainoperad,C^{\otimes})$, the map
$y=\left(\prod_{\mathbf{u}\in\mathrf{P}(k)}y(\mathbf{u})\right)\circ\sethom(1_{\mainoperad},\kappa^{\otimes})$,
where \begin{multline*}
y(\mathbf{u})=\bar{y}(\mathbf{u})|\sethomn{\mathbf{u}}(\mathcal{V}(\mathbf{u}),P^{\mathbf{u}}):\sethomn{\mathbf{u}}(\mathcal{V}(\mathbf{u}),P^{\mathbf{u}})\\
\to\homz(\mathcal{V}(\mathbf{u})\otimes\bigotimes_{\mathbf{u}}\mathcal{V}(u_{i}),C^{\mathfrak{g}(\mathbf{u})})\end{multline*}
 and the maps\[
\bar{y}(\mathbf{u}):\homz(\mathcal{V}(\mathbf{u}),P^{\mathbf{u}})\to\homz(\mathcal{V}(\mathbf{u})\otimes\bigotimes_{\mathbf{u}}\mathcal{V}(u_{i}),C^{\mathfrak{g}(\mathbf{u})})\]
 map the factor \[
\homz(\mathcal{V}(\mathbf{u}),\bigotimes_{\mathbf{u}}L(u_{j}))\subset\homz(\mathcal{V}(\mathbf{u}),P^{\mathbf{u}})\]
with $L(u_{j})=\homz(\mathcal{V}(u_{j}),C^{u_{j}})$ via the map induced
by the associativity of the $\mathrm{Hom}$ and $\otimes$ functors.
\end{enumerate}
Consider the diagram whose rows are copies of diagram~\ref{eq:bigdiag}
\begingroup\small\begin{equation}\xymatrix@C-5pt{{K}\ar@{^{(}->}[d]\ar[r]^{g|K} & {W}\ar@{=}[d] & {\sethom(\mainoperad,(L_{\mainoperad}C)^\otimes)}\ar@{_{(}->}[l]_-{y} \ar@{=}[d]\\ {L_{\mainoperad}C}\ar[r]^-{g}\ar[d]_{\kappa} & {W}\ar[d]^{r_1} & {\sethom(\mainoperad,(L_{\mainoperad}C)^\otimes)}\ar@{_{(}->}[l]_-{y}\ar[d]^{\sethom(\iota,1)}  \\ {\sethom(\mathcal{I},C^\otimes)}\ar[r]_-{g_1} & {T} & {\sethom(\mathcal{I},(L_{\mainoperad}C)^\otimes)}\ar@{_{(}->}[l]^-{w}}\label{eq:3dia}\end{equation}\endgroup
where

\begin{enumerate}
\item $W={\displaystyle \prod_{\mathbf{u}\in\mathrf{P}(k)}\homz(\mathcal{H}(\mathbf{u})\otimes\bigotimes_{u_{i}\in\mathbf{u}}\mathcal{H}(u_{i}),C^{\mathfrak{g}(\mathbf{u})})}$
\item $y:\sethom(\mainoperad,(L_{\mainoperad}C)^{\otimes})\to W$ is defined
as in diagram~\ref{eq:bigdiag}.
\item $T={\displaystyle \prod_{\mathbf{u}\in\mathrf{P}(k)}\homz(\mathcal{I}(\mathbf{u})\otimes\bigotimes_{u_{i}\in\mathbf{u}}\mathcal{H}(u_{i}),C^{\mathfrak{g}(\mathbf{u})})}$
\item The map $g_{1}=\left(\prod_{\mathbf{u}\in\mathrf{P}(k)}\homz(\bar{\gamma}(\mathbf{u}),1)\right)\circ\kappa$
where \[
\bar{\gamma}(\mathbf{u})=\gamma(\mathbf{u})|{\displaystyle \mathcal{I}(\mathbf{u})\otimes\bigotimes_{u_{i}\in\mathbf{u}}\mainoperad(u_{i})}:{\displaystyle \mathcal{I}(\mathbf{u})\otimes\bigotimes_{u_{i}\in\mathbf{u}}\mainoperad(u_{i})}\to\mainoperad(\mathfrak{g}(\mathbf{u}))\]

\item The map $r_{1}=\homz(j_{1},1)$ where\[
j_{1}:\mathcal{I}(\mathbf{u})\otimes\bigotimes_{u_{i}\in\mathbf{u}}\mainoperad(u_{i})\hookrightarrow\mathcal{H}(\mathbf{u})\otimes\bigotimes_{u_{i}\in\mathbf{u}}\mainoperad(u_{i})\]
for $\mathbf{u}\in\mathrf{P}(k)$, are the inclusions.
\item $\iota:\mathcal{I}\hookrightarrow\mathcal{H}$ is the inclusion.
\end{enumerate}
Suppose $r\in K$. Then the image of $r$ under the downward maps
on the left of diagram~\ref{eq:3dia} must be $0$, since $K$ is
the kernel of $\kappa$. On the other hand, $r=y(h)$ in the top two
rows of this diagram.

The commutativity of diagram~\ref{eq:3dia} implies that the image
of $h$ under the downward maps on the \emph{right} is also $0$,
so that the coproduct of $r$ in the kernel of $\sethom(\iota,1)$.
This implies that the coproduct of $K$ is the pullback of a map \[
K\to\sethom(\mathcal{H}/\mathcal{I},C^{\otimes})\]
over the projection $p:\mainoperad\to\mainoperad/\mathcal{I}$. 

Let $X\in\setf$ and let $x\in X$ be an arbitrary element. We claim
that the diagram

\begingroup\small\begin{equation}\xymatrix@C+10pt{{K}\ar@{^{(}->}[d]\ar[r]^{g|K} & {W'}\ar@{^{(}->}[d]^{\homz(p\otimes1,1)} & {\sethom(\mainoperad/\mathcal{I},(L_{\mainoperad}C)^\otimes)}\ar@{_{(}->}[l]_-{y} \ar@{^{(}->}[d]^{\sethom(p,1)}\\ {L_{\mainoperad}C}\ar[r]^-{g}\ar[d]_{\kappa} & {W}\ar[d]^{ p_{W(X)}} & {\sethom(\mainoperad,(L_{\mainoperad}C)^\otimes)}\ar@{_{(}->}[l]_-{y}\ar[d]^{p_X}  \\ {\sethom(\mainoperad,C^\otimes)}\ar[d]_{\homz(\iota,1)}\ar[r]^-{\homz(\gamma,1)}& {W(X)}\ar[d]^{\theta(X,x)}&{\sethomn{X}(\mainoperad(X),(L_{\mainoperad}(C))^X)}\ar@{_{(}->}[l]_-{y(X)}\ar[d]^{\phi(X,x)}\\ {\sethom(\mathcal{I},C^\otimes)}\ar[r]_-{\homz(\gamma,1)} & { Y(X,x)} & { \sethomn{X}(\mainoperad(X),M(X,x))}\ar@{_{(}->}[l]^-{ y(X,x)}}\label{eq:4dia}\end{equation}\endgroup 

commutes, where

\begin{enumerate}
\item $p_{X}$ and $p_{W(X)}$ are projections onto direct factors.
\item $W'={\displaystyle \prod_{\mathbf{u}\in\mathrf{P}(k)}\homz(\mathcal{H}(\mathbf{u})/\mathcal{I}(\mathbf{u})\otimes\bigotimes_{u_{i}\in\mathbf{u}}\mathcal{H}(u_{i}),C^{\mathfrak{g}(\mathbf{u})})}$
\item $W(X)={\displaystyle \prod_{\substack{\mathbf{u}\in\mathrf{P}(k)\\
s:X\to\mathfrak{f}(\mathbf{u})}
}}\homz(\mainoperad(\mathbf{u})\otimes{\displaystyle \bigotimes_{u_{i}\in\mathbf{u}}\mainoperad(u_{i}),C^{\mathfrak{g}(\mathbf{u})})}$, where $s:X\to\mathfrak{f}(\mathbf{u})$ is a set-bijection. This
is exactly like $W$, except that we only consider $\mathbf{u}$ such
that $\mathfrak{f}(\mathbf{u})$ has the same \emph{cardinality} as
the set $X$.
\item ${\displaystyle Y(X,x)=\prod_{\substack{\mathbf{u}\in\mathrf{P}(k)\\
s:X\to\mathfrak{f}(\mathbf{u})}
}\homz(\mainoperad(\mathbf{u})\otimes\bigotimes_{\mathbf{u},s(x)}(\mainoperad,\mathcal{I}),C^{\mathfrak{g}(\mathbf{u})})}$, where $x\in X$ is any element --- see definition~\ref{def:unorderedspans}.
This is exactly like $W(X)$, except that the $s(x)^{\text{th}}$
factor of $\mainoperad(u_{i})$ has been replaced with $\mathcal{I}(s(x))$.
\item $\theta(X,x)={\displaystyle \homz(1\otimes\bigotimes_{X,x}(1,\iota),1):W(X)\to Y(X,x)}$.
This is the dual of $1\otimes\bigotimes_{X,x}(1,\iota)$, which is
the identity, except for the $x^{\text{th}}$ factor on the right.
For this factor it is the inclusion $\iota:\mathcal{I}(s(x))\hookrightarrow\mainoperad(s(x))$.
\item $\phi(X,x)=\homz\left(1,{\displaystyle \bigotimes_{X,x}(1,\homz(\iota,1))}\right)$
\item ${\displaystyle M(X,x)=\bigotimes_{X,x}(\sethom(\mainoperad,C^{\otimes}),\sethom(\mathcal{I},C^{\otimes}))}$
--- see definition~\ref{def:unorderedspans}. This is the similar
to $(\sethom(\mainoperad,C^{\otimes})^{X}$, except that the $x^{\text{th}}$
factor has been replaced by $\sethom(\mathcal{I},C^{\otimes})$.
\item $y(X)$ is defined as $y(\mathbf{u})$ in diagram~\ref{eq:bigdiag}
and $y(X,x)$ is defined analogously --- with the $x^{\text{th}}$
factor mapping $\sethom(\mathcal{I},C^{\otimes})$.
\end{enumerate}
The upper squares of diagram~\ref{eq:4dia} commute because they
did in diagram~\ref{eq:3dia}.

The lower left square of diagram~\ref{eq:4dia} commutes because
it is the dual of the diagram \[\xymatrix{{\mainoperad(\mathfrak{g}(\mathbf{u}))} & {\mainoperad(\mathbf{u})\otimes\bigotimes_{\mathbf{u}}\mainoperad(u_i)}\ar[l]_-{\gamma} \\ {\mathcal{I}(\mathfrak{g}(\mathbf{u}))}\ar[u]^{\iota} & {\mainoperad(\mathbf{u})\otimes\bigotimes_{\mathbf{u},s(x)}(\mainoperad,\mathcal{I})}\ar[l]^-{\gamma}\ar[u]_{1\otimes\bigotimes_{\mathbf{u},s(x)}(1,\iota)}}\]
which is well-defined and commutes because $\mathcal{I}$ is an operadic
ideal of $\mainoperad$. The lower right square of diagram~\ref{eq:4dia}
commutes because of the naturality of the $y$-maps.

A diagram-chase around the outer rim of diagram~\ref{eq:4dia} shows
that if $k\in K$, then the coproduct of $k$, evaluated on any element
of $\mathcal{H}(X)/\mathcal{I}(X)$ (or $\mathcal{H}(X)$) gives a
result that lies in the kernel of ${\displaystyle \bigotimes_{X,x}(1,\homz(\iota,1))}$
for any finite set $X$ and any element $x\in X$, hence is in $K^{X}$
--- see proposition~\ref{pro:kerunorderedspan}. It follows that
$K$ is a sub-coalgebra of $L_{\mainoperad}C$ and one that has been
pulled back from $\mainoperad/\mathcal{I}$.

The lemma's final statement follows from the universal property of
cofree coalgebras. Suppose $M$ is any coalgebra over $\mainoperad/\mathcal{I}$
equipped with a chain-map $\alpha:M\to C$. By composition with the
projection $p:\mainoperad\to\mainoperad/\mathcal{I}$, we may regard
$M$ as a coalgebra over $\mainoperad$. The universal property of
a cofree coalgebra implies that there exists a \emph{unique} morphism
of $\mainoperad$-coalgebras\[
M\to L_{\mainoperad}C\]
 that makes the diagram \[\xymatrix{{M}\ar[r]\ar[rd]_{\alpha} & {L_{\mainoperad}C}\ar[d]^{\varepsilon} \\ {} & {C}}\]
commute (where $\varepsilon:L_{\mainoperad}C\to C$ is the cogeneration
map). But the image of $M$ must lie in $K\subseteq L_{\mainoperad}C$,
hence $K$ has the universal property of a cofree coalgebra over $\mainoperad/\mathcal{I}$.{

}

%%bibliography begins
\providecommand{\bysame}{\leavevmode\hbox to3em{\hrulefill}\thinspace}
\providecommand{\MR}{\relax\ifhmode\unskip\space\fi MR }
\providecommand{\MRhref}[2]{%
  \href{http://www.ams.org/mathscinet-getitem?mr=#1}{#2}
}
\providecommand{\href}[2]{#2}

%% index begins


\begin{thebibliography}{10}

\bibitem{Baer:1937}
R.~Baer, \emph{Abelian groups without elements of finite order}, Duke Math. J.
  \textbf{3} (1937), 68--122.

\bibitem{Berger-moerdijk-axiom-operad}
Clemens Berger and Ieke Moerdijk, \emph{Axiomatic homotopy theory for operads},
  Comment. Math. Helv. \textbf{78} (2003), no.~4, 681--721.

\bibitem{Blass-Gobel:1996}
Andreas~R. Blass and R{\"u}diger G{\"o}bel, \emph{Subgroups of the
  {B}aer-{S}pecker group with few endomorphisms but large dual}, Fundamenta
  Mathematicae \textbf{149} (1996), 19--29.

\bibitem{Baer-Specker-survey}
Eoin Coleman, \emph{The {B}aer-{S}pecker group}, Bulletin of the Irish
  Mathematical Society (1998), no.~40, 9--23.

\bibitem{Fuchs:1970}
L.~Fuchs, \emph{Abelian groups}, vol. {I} and {II}, Academic Press, 1970 and
  1973.

\bibitem{Gugenheim:1960}
V.~K. A.~M. Gugenheim, \emph{On a theorem of {E.} {H.} {B}rown}, Illinois J. of
  Math. \textbf{4} (1960), 292--311.

\bibitem{Kriz-May}
I.~Kriz and J.~P. May, \emph{Operads, algebras, modules and motives},
  Ast{\'e}risque, vol. 233, Soci{\'e}t{\'e} {M}ath{\'e}matique de {France},
  1995.

\bibitem{MacLane:cw}
S.~MacLane, \emph{Categories for the {W}orking {M}athematician},
  Springer-Verlag, Berlin, Heidelberg, New York, 1971.

\bibitem{operad-book}
Martin Markl, Steve Shnider, and Jim Stasheff, \emph{Operads in {A}lgebra,
  {T}opology and {P}hysics}, Mathematical Surveys and Monographs, vol.~96,
  American Mathematical Society, May 2002.

\bibitem{Smith:1994}
Justin~R. Smith, \emph{Iterating the cobar construction}, vol. 109, Memoirs of
  the A. M. S., no. 524, American Mathematical Society, Providence, Rhode
  Island, May 1994.

\bibitem{Smith-cofree}
\bysame, \emph{Cofree coalgebras over operads}, Topology and its Applications
  \textbf{133} (2003), 105--138.

\bibitem{Baer-Specker-nonfree}
E.~Specker, \emph{Additive {G}ruppen von {F}olgen ganzer {Z}ahlen}, Portugaliae
  Math. \textbf{9} (1950), 131--140.

\bibitem{Sweedler:1969}
Moss~E. Sweedler, \emph{Hopf algebras}, W. A. Benjamin, Inc., 1969.

\bibitem{Weibel:homological-algebra}
Charles~A. Weibel, \emph{An introduction to homological algebra}, Cambridge
  studies in advanced mathematics, vol.~38, Cambridge University Press, 1994.

\end{thebibliography}
    \end{document}